# Fourier series (based) multiscale method for computational analysis in science and engineering:

# IV. Fourier series multiscale solution for the convection-diffusion-reaction equation

Weiming Sun[a*+] and Zimao Zhang[b]

**Abstract:** Fourier series multiscale method, a concise and efficient analytical approach for multiscale computation, will be developed out of this series of papers. In the fourth paper, the application of the Fourier series multiscale method to the one- and two-dimensional convection-diffusion-reaction equations is implemented, where the Fourier series multiscale solutions are derived, the convergence characteristics of the Fourier series multiscale solutions are investigated by numerical examples, the schemes for application of the Fourier series multiscale method are optimized, and the multiscale characteristics of the convection-diffusion-reaction equations are demonstrated. The preliminary study on applications verifies the effectiveness of the present Fourier series multiscale method and provides a reliable reference which can be used for persistent improvement in computational performance of other multiscale methods.

## 1. Introduction

With regard to the researches on multiscale methods of computational analysis in science and engineering, the convection-diffusion-reaction equation plays an important role in theory and practice. On the one hand, the convection-diffusion-reaction equation and its particular cases (e.g. convection-diffusion equation, Helmholtz equation, convection-reaction equation and diffusion-reaction equation) are widely applied to the simulations of many physical and chemical processes [1, 2] such as turbulence, chemistry, combustion, heat transfer with radiation, diffusion of pollutants and so on, which profoundly reveals the diversity of the

[a]Department of Mathematics and Big Data, School of Artificial Intelligence, Jianghan University, Wuhan, 430056, China

[b]Department of Mechanics, School of Civil Engineering, Beijing Jiaotong University, Beijing, 100044, China

*Correspondence to: Weiming Sun, Department of Mathematics and Big Data, School of Artificial Intelligence, Jianghan University, Wuhan, 430056, China

[+] E-mail: xuxinenglish@hust.edu.cn

scientific discipline background of the multiscale phenomena and the multiscale problems and presents a fresh point of view for the deeper understanding of the complexity of natural phenomena and engineering systems. On the other hand, the convection-diffusion-reaction equation has varied physical behaviors with the presence of the exponential regime and the propagation regime [1-3], of which the solving process characterizes several essential contradictions between the multiscale problems and the multiscale methods, such as selection of computational scale and computational feasibility, tuning of computational parameters and stability of computational performance, capturing of localized multiscale structures and computational accuracy, and also points out the research directions in developing the practical, feasible, efficient and highly accurate multiscale methods.

Nowadays, in order to meet the urgent needs of solving the convection-diffusion-reaction equation, the researchers have developed a number of multiscale methods, such as stabilized finite element method, bubble method, wavelet based finite element method, meshless method, finite increment calculus based multiscale method, and variational multiscale method, which bring to dramatic advancements for the current multiscale analysis methods [4]. All these multiscale methods have in common that they are based on the traditional numerical methods (especially the Galerkin finite element method) and developed by proper modification of the original computational procedures or the routine application schemes. For example, the stabilized finite element method is the first computational method proposed for solving the multiscale problem of the convection-diffusion equation. As is well known, the idea of constructing this kind of multiscale method is very simple, which corresponds to introducing additional stabilization terms to the Galerkin variational formulation. The additional stabilization terms are the weighted residuals within each element of the differential equation being solved. Accordingly, different variations on the form of weighting functions lead to different improvements on stability behavior, and therefore, different versions of the stabilized finite element methods [5-11]. Instead of a modification of the Garlekin variational formulation, the bubble function method is developed with the finite element space enriched by special functions, denoted by bubbles for simplicity. In this process, the design of bubble functions becomes a key to the construction of the multiscale method and several bubble designs, including the classic bubble function [12-14], the residual-free bubble function [15-22] and the multiscale function [23], have been successively suggested. Similarly, by enriching the finite element space with the wavelets [24] as basis functions, the researchers construct another effective multiscale method, namely, the wavelet based finite element method [25] for the analysis of the convection-diffusion-reaction equation. The wavelet based finite element method is capable of adjusting adaptively the scales of analysis and validated with enhanced numerical stability, high computational speed and high accuracy [4, 26, 27]. The meshless method, also known as element free method, is a new type numerical method developed widely in recent years [28, 29]. This method yields the localized, accurate meshless approximation from a set of nodal data and associated weight functions with compact support on the domain. Therefore, it bypasses the significant difficulties which the use of the finite element method presents for the analysis of sharp gradient and singularity zones, and possesses the necessary intrinsic adaptability to multiscale problems such as the convection-diffusion-reaction equation. The finite calculus (FIC) method newly developed by Onate and co-workers is a consistent procedure to re-formulate the classic differential equation governing the multiscale phenomena [30-32]. The derived FIC formulation introduces naturally new terms involving characteristic length parameters. Application of the Galerkin finite element method to the FIC formulation leads to a novel, conceptually different multiscale method [3, 33-35]. By contrast, the variational multiscale method (VMS), as originally proposed by Hughes in 1995, offers a fresh theoretical framework for the solution of multiscale problems [36]. The implementation of this framework consists of two steps: the construction of the subgrid scale model and its solution with the Galerkin finite element method. The first step is non-numerical and involves a series of tasks, such as introduction of the scale groups with coarse scale and fine (subgrid) scale, decomposition of the solution into coarse and fine scale components, solution of the fine scale equations in term of the coarse scale residual, and elimination of the fine scale solution from the coarse scale equation. This procedure leads to a modified coarse scale equation, sometimes referred to a subgrid scale model. The resulting subgrid scale model turns out to be a suitable one for presentation to the Galerkin finite element method and accurate numerical results naturally follow in the second step of the variational multiscale method. Therefore, the variational multiscale method makes the process of solving multiscale problems concise and consistent again [37, 38].

Over the years, these multiscale methods have found wide applications to a variety of multiscale problems, for instance, the typical convection-diffusion-reaction equation. However, it is acknowledged that the current researches on multiscale methods are still limited within the theoretical framework of numerical methods, and inevitably, these developed multiscale methods have disadvantages such as high costs of computation, low accuracy for higher order derivatives of field variables, and difficulties in analysis of computational parameters' effects on computational results [28].

Therefore, in the fourth paper of the series of researches on Fourier series multiscale method, we have taken a strategic shift from the theoretical framework of numerical methods towards that of the Fourier series multiscale method and performed a reinvestigation of the convection-diffusion-reaction equation. In this paper, the Fourier series multiscale solutions of the one-dimensional and two-dimensional convection-diffusion-reaction equations are derived. Next, convergence characteristics of the Fourier series multiscale solutions are investigated by numerical examples. With settings of computational schemes of the Fourier series multiscale method optimized, the multiscale characteristics of the one-dimensional and two-dimensional convection-diffusion-reaction equations are finally demonstrated.

## 2. Fourier series multiscale solution for one-dimensional convection-diffusion-reaction equation

The one-dimensional convection-diffusion-reaction equation is a second order ($2r = 2$) linear differential equation with constant coefficients, where convection refers to the term involving first order derivative of the unknown function. Therefore, in this section we rationally specify the type of the Fourier series expansion of the one-dimensional solution as the full-range Fourier series and closely follow the step outlined in [39] to derive the Fourier series multiscale solution of the equation.

*2.1. Description of the problem*

With reference to [2, 20, 21], the one-dimensional convection-diffusion-reaction equation on the interval $[-a, a]$ can be written in the dimensionless form as

$$\mathcal{L}_{cdr,1} \varphi = f , \tag{1}$$

where the differential operator

$$\mathcal{L}_{cdr,1} = P_e \frac{d}{dx_1} - \frac{d^2}{dx_1^2} - P_e D_a , \tag{2}$$

$\varphi(x_1)$ is the unknown function, $f(x_1)$ is a given source function, $P_e$ is the dimensionless

Peclet number related to the convective velocity, and $D_a$ is the dimensionless Damkohler number.

The boundary $x_1 = a$ (actually an endpoint of the solution interval $[-a, a]$) is taken as an example and the corresponding boundary conditions are classified into the following two types:

1. Dirichlet boundary condition (for short, D boundary condition)
$$\varphi(a) = \bar{\varphi}_a, \tag{3}$$
2. Neumann boundary condition (for short, N boundary condition)
$$\varphi^{(1)}(a) = \bar{\phi}_a, \tag{4}$$

where $\bar{\varphi}_a$ and $\bar{\phi}_a$ are the prescribed values respectively of the unknown function and its first order derivative at the boundary $x_1 = a$.

## 2.2. The general solution

Suppose that Eq. (1) has homogeneous solutions of the following form
$$p_H(x_1) = \exp(\eta x_1), \tag{5}$$
where $\eta$ is an undetermined constant.

Substituting it into the homogeneous form of Eq. (1), we obtain the characteristic equation
$$\eta^2 - P_e \eta + P_e D_a = 0. \tag{6}$$

It is obvious that Eq. (6) has two distinct real roots when $4D_a < P_e$, two distinct complex roots when $4D_a > P_e$ and one double real root when $4D_a = P_e$.

Then we can further express the characteristic roots as
$$\eta_1 = -\alpha_{10} + \alpha_{30} - i\alpha_{20} \quad \text{and} \quad \eta_2 = -\alpha_{10} - \alpha_{30} + i\alpha_{20}, \tag{7}$$
where $i = \sqrt{-1}$.

Table 1 presents a detailed distribution of real constants $\alpha_{10}$, $\alpha_{20}$ and $\alpha_{30}$.

Accordingly, Table 2 presents the two linearly independent homogeneous solutions $\{p_{l,H}(x_1)\}_{1 \leq l \leq 2}$ of Eq. (1).

Table 1: Distribution of real constants in the characteristic roots.

|  | $\alpha_{10}$ | $\alpha_{20}$ | $\alpha_{30}$ |
| --- | --- | --- | --- |
| $4D_a < P_e$ | $-\dfrac{P_e}{2}$ | 0 | $\dfrac{\sqrt{P_e^2 - 4P_e D_a}}{2}$ |
| $4D_a > P_e$ | $-\dfrac{P_e}{2}$ | $\dfrac{\sqrt{4P_e D_a - P_e^2}}{2}$ | 0 |
| $4D_a = P_e$ | $-\dfrac{P_e}{2}$ | 0 | 0 |

Table 2: Expressions for $\{p_{l,H}(x_1)\}_{1\leq l\leq 2}$.

|  | $p_{1,H}(x_1)$ | $p_{2,H}(x_1)$ |
|---|---|---|
| $4D_a < P_e$ | $\exp[(-\alpha_{10}+\alpha_{30})x_1]$ | $\exp[(-\alpha_{10}-\alpha_{30})x_1]$ |
| $4D_a > P_e$ | $\exp(-\alpha_{10}x_1)\sin(\alpha_{20}x_1)$ | $\exp(-\alpha_{10}x_1)\sin[\alpha_{20}(a-x_1)]$ |
| $4D_a = P_e$ | $\exp(-\alpha_{10}x_1)$ | $x_1\exp(-\alpha_{10}x_1)$ |

We select the following vector of functions
$$\mathbf{p}_1^T(x_1) = [p_{1,H}(x_1) \quad p_{2,H}(x_1)], \tag{8}$$
and then we can express the general solution of Eq. (1) as
$$\varphi_1(x_1) = \mathbf{\Phi}_1^T(x_1) \cdot \mathbf{q}_1, \tag{9}$$
where the definitions of the vector of basis function $\mathbf{\Phi}_1^T(x_1)$ and the vector of undetermined constants $\mathbf{q}_1$ are respectively referred to Eqs. (22) and (20) in [40].

*2.3. The supplementary solution*

Let $N_1^s$ be a positive integer and $\Upsilon^s = \{x_{1,n_1}, n_1 = 1, 2, \cdots, N_1^s + 1\}$ be a set of interpolation points uniformly distributed on the one-dimensional interval $[-a, a]$. Now we are to construct an interpolation algebraical polynomial $f_s(x_1)$ of the source function $f(x_1)$, such that the following interpolation conditions are satisfied
$$f_s(x_{1,n_1}) = f(x_{1,n_1}), \quad n_1 = 1, 2, \cdots, N_1^s + 1. \tag{10}$$
For this purpose, we select the vector of functions in the form of algebraical polynomials
$$\mathbf{p}_{fs}^T(x_1) = [p_{fs,1}(x_1) \quad p_{fs,2}(x_1) \quad \cdots \quad p_{fs,N_1^s+1}(x_1)], \tag{11}$$
where
$$p_{fs,j}(x_1) = (x_1/a)^{j-1}, \quad j = 1, 2, \cdots, N_1^s + 1, \tag{12}$$
and define the vector of undetermined constants
$$\mathbf{a}_{fs}^T = [H_{fs,1} \quad H_{fs,2} \quad \cdots \quad H_{fs,N_1^s+1}], \tag{13}$$
then the interpolation algebraical polynomial $f_s(x_1)$ can be expressed as
$$f_s(x_1) = \mathbf{p}_{fs}^T(x_1) \cdot \mathbf{a}_{fs}. \tag{14}$$
Combining the interpolation conditions in Eq. (10), we obtain
$$\mathbf{R}_{fs}\mathbf{a}_{fs} = \mathbf{q}_{fs}, \tag{15}$$
where
$$\mathbf{R}_{fs} = \begin{bmatrix} \mathbf{p}_{fs}^T(x_{1,1}) \\ \mathbf{p}_{fs}^T(x_{1,2}) \\ \vdots \\ \mathbf{p}_{fs}^T(x_{1,N_1^s+1}) \end{bmatrix}, \tag{16}$$
and the vector of values of the source function
$$\mathbf{q}_{fs}^T = [f(x_{1,1}) \quad f(x_{1,2}) \quad \cdots \quad f(x_{1,N_1^s+1})]. \tag{17}$$

Further, we consider the supplementary solution of Eq. (1) in the form of algebraical polynomials.

1. For $P_e \neq 0$ and $D_a \neq 0$, we select, without loss of generality, the vector of functions in the form of algebraical polynomials

$$\mathbf{p}_s^T(x_1) = [p_{s,1}(x_1) \quad p_{s,2}(x_1) \quad \cdots \quad p_{s,N_1^s+1}(x_1)], \tag{18}$$

where

$$p_{s,j}(x_1) = (x_1/a)^{j-1}, \quad j = 1, 2, \cdots, N_1^s + 1, \tag{19}$$

and define the vector of undetermined constants

$$\mathbf{a}_s^T = [G_{s,1} \quad G_{s,2} \quad \cdots \quad G_{s,N_1^s+1}], \tag{20}$$

thus, the supplementary solution of Eq. (1) can be expressed as

$$\varphi_s(x_1) = \mathbf{p}_s^T(x_1) \cdot \mathbf{a}_s. \tag{21}$$

Substituting the supplementary solution into the equation

$$\mathcal{L}_{cdr,1}\varphi_s = f_s, \tag{22}$$

we obtain

$$\mathbf{R}_s \mathbf{a}_s = \mathbf{a}_{fs}, \tag{23}$$

where the transformation matrix

$$\mathbf{R}_s = \begin{bmatrix} -P_e D_a & P_e/a & -1 \cdot 2/a^2 & 0 & \cdots & 0 \\ 0 & -P_e D_a & 2P_e/a & -2 \cdot 3/a^2 & \cdots & 0 \\ \vdots & \vdots & \ddots & \ddots & \ddots & \vdots \\ 0 & 0 & 0 & -P_e D_a & (N_1^s - 1)P_e/a & -(N_1^s - 1)N_{x_1}/a^2 \\ 0 & 0 & 0 & \cdots & -P_e D_a & N_1^s P_e/a \\ 0 & 0 & 0 & \cdots & 0 & -P_e D_a \end{bmatrix}. \tag{24}$$

Hence

$$\mathbf{a}_s = \mathbf{R}_s^{-1} \mathbf{a}_{fs}. \tag{25}$$

Define the vector of basis functions

$$\boldsymbol{\Phi}_s^T(x_1) = \mathbf{p}_s^T(x_1) \cdot \mathbf{R}_s^{-1} \mathbf{R}_{fs}^{-1}. \tag{26}$$

Then the supplementary solution $\varphi_s(x_1)$ can be expressed as

$$\varphi_s(x_1) = \boldsymbol{\Phi}_s^T(x_1) \cdot \mathbf{q}_{fs}. \tag{27}$$

2. For $P_e \neq 0$ and $D_a = 0$, we replace the components in $\mathbf{p}_s^T(x_1)$ with

$$p_{s,j}(x_1) = (x_1/a)^j, \quad j = 1, 2, \cdots, N_1^s + 1, \tag{28}$$

and the transformation matrix with

$$\mathbf{R}_s = \begin{bmatrix} P_e/a & -1 \cdot 2/a^2 & 0 & \cdots & 0 \\ 0 & 2P_e/a & -2 \cdot 3/a^2 & \cdots & 0 \\ \vdots & \vdots & \ddots & \ddots & \vdots \\ 0 & 0 & \cdots & N_1^s P_e/a & -N_1^s(N_1^s+1)/a^2 \\ 0 & 0 & \cdots & 0 & (N_1^s+1)P_e/a \end{bmatrix}, \tag{29}$$

then the supplementary solution $\varphi_s(x_1)$ can be expressed by Eq. (27).

3. For $P_e = 0$ and $D_a = 0$, we replace the components in $\mathbf{p}_s^T(x_1)$ with
$$p_{s,j}(x_1) = (x_1/a)^{j+1}, \quad j = 1, 2, \cdots, N_1^s + 1, \tag{30}$$
and the transformation matrix with
$$\mathbf{R}_s = \begin{bmatrix} -1 \cdot 2/a^2 & 0 & \cdots & 0 \\ 0 & -2 \cdot 3/a^2 & \cdots & 0 \\ \vdots & \vdots & \ddots & \vdots \\ 0 & 0 & 0 & -(N_1^s+1)(N_1^s+2)/a^2 \end{bmatrix}, \tag{31}$$
then the supplementary solution $\varphi_s(x_1)$ can be expressed by Eq. (27).

In some occasions the supplementary solution $\varphi_s(x_1)$ is not introduced into the Fourier series multiscale solution, $N_1^s = 0$ is denoted for convenience.

## 2.4. The particular solution

Let the error of the interpolation algebraical polynomial $f_s(x_1)$ relative to the source function $f(x_1)$ be
$$f_p(x_1) = f(x_1) - f_s(x_1), \tag{32}$$
then with the interpolation conditions we observe that
$$f_p(-a) = 0, \quad f_p(a) = 0. \tag{33}$$
Expand it in a full-range Fourier series on the interval $[-a, a]$, we have
$$f_p(x_1) = \sum_{m=0}^{\infty} \mu_m [V_{fp,1m} \cos(\alpha_m x_1) + V_{fp,2m} \sin(\alpha_m x_1)], \tag{34}$$
where $\alpha_m = m\pi/a$, $\mu_m = \begin{cases} 1/2, & m = 0 \\ 1, & m > 0 \end{cases}$, and $V_{fp,1m}$, $V_{fp,2m}$ are the Fourier coefficients of $f_p(x_1)$.

Further, $f_p(x_1)$ can be expressed in matrix form
$$f_p(x_1) = \mathbf{\Phi}_0^T(x_1) \cdot \mathbf{q}_{fp}, \tag{35}$$
where the definitions of the vector of one-dimensional trigonometric functions $\mathbf{\Phi}_0^T(x_1)$ and the vector of Fourier coefficients $\mathbf{q}_{fp}$ are respectively referred to Eqs. (33) and (39) in [40].

Meanwhile, suppose that Eq. (1) has the particular solution in the form of Fourier series
$$\varphi_0(x_1) = \mathbf{\Phi}_0^T(x_1) \cdot \mathbf{q}_0, \tag{36}$$
where $\mathbf{q}_0$ is the vector of undetermined Fourier coefficients.

Substitute Eqs. (35) and (36) into the equation
$$\mathcal{L}_{cdr,1} \varphi_0 = f_p, \tag{37}$$
then we obtain
$$[\mathcal{L}_{cdr,1} \mathbf{\Phi}_0^T(x_1)] \cdot \mathbf{q}_0 = \mathbf{\Phi}_0^T(x_1) \cdot \mathbf{q}_{fp}, \tag{38}$$
which can be solved by the collocation method (CM) or the Fourier coefficient comparison method (FCCM).

1. The collocation method

For the particular solution $\varphi_0(x_1)$, let $M$ be the number of truncated terms of the Fourier series, and we denote the set of collocation points uniformly distributed on the one-dimensional interval $[-a, a]$ by $\Upsilon_p = \{x_1^{n_1}, n_1 = 1, 2, \cdots, 2M+1\}$. Substituting these collocation points into Eq. (38), we obtain

$$\mathbf{R}_{p1}\mathbf{q}_0 = \mathbf{R}_{p2}\mathbf{q}_{fp}, \tag{39}$$

where

$$\mathbf{R}_{p1} = \begin{bmatrix} P_e \mathbf{\Phi}_0^{(1)T}(x_1^1) - \mathbf{\Phi}_0^{(2)T}(x_1^1) - P_e D_a \mathbf{\Phi}_0^T(x_1^1) \\ P_e \mathbf{\Phi}_0^{(1)T}(x_1^2) - \mathbf{\Phi}_0^{(2)T}(x_1^2) - P_e D_a \mathbf{\Phi}_0^T(x_1^2) \\ \vdots \\ P_e \mathbf{\Phi}_0^{(1)T}(x_1^{2M+1}) - \mathbf{\Phi}_0^{(2)T}(x_1^{2M+1}) - P_e D_a \mathbf{\Phi}_0^T(x_1^{2M+1}) \end{bmatrix}, \tag{40}$$

$$\mathbf{R}_{p2} = \begin{bmatrix} \mathbf{\Phi}_0^T(x_1^1) \\ \mathbf{\Phi}_0^T(x_1^2) \\ \vdots \\ \mathbf{\Phi}_0^T(x_1^{2M+1}) \end{bmatrix}, \tag{41}$$

and therefore

$$\mathbf{q}_0 = \mathbf{R}_{p1}^{-1} \mathbf{R}_{p2} \mathbf{q}_{fp}. \tag{42}$$

2. The Fourier coefficient comparison method

For the particular solution $\varphi_0(x_1)$, let $M$ be the number of truncated terms of the Fourier series, and compare the first $2M+1$ Fourier coefficients on the both sides of Eq. (38) successively, then

a. for $m = 0$, we have

$$-P_e D_a V_{10} = V_{fp,10}, \tag{43}$$

and therefore

$$V_{10} = -\frac{1}{P_e D_a} \cdot V_{fp,10}. \tag{44}$$

b. For $m = 1, 2, \cdots, M$, we have

$$\begin{bmatrix} \alpha_m^2 - P_e D_a & P_e \alpha_m \\ -P_e \alpha_m & \alpha_m^2 - P_e D_a \end{bmatrix} \begin{bmatrix} V_{1m} \\ V_{2m} \end{bmatrix} = \begin{bmatrix} V_{fp,1m} \\ V_{fp,2m} \end{bmatrix}, \tag{45}$$

and therefore

$$\begin{bmatrix} V_{1m} \\ V_{2m} \end{bmatrix} = \begin{bmatrix} \alpha_m^2 - P_e D_a & P_e \alpha_m \\ -P_e \alpha_m & \alpha_m^2 - P_e D_a \end{bmatrix}^{-1} \begin{bmatrix} V_{fp,1m} \\ V_{fp,2m} \end{bmatrix}. \tag{46}$$

By combining Eqs. (44) and (46), we obtain the vector of undetermined Fourier coefficients $\mathbf{q}_0$.

## 2.5. The Fourier series multiscale solution

Putting the Eqs. (36), (9) and (27) together, we thus express the Fourier series multiscale solution of the one-dimensional convection-diffusion-reaction equation as

$$\begin{aligned} \varphi(x_1) &= \varphi_0(x_1) + \varphi_1(x_1) + \varphi_s(x_1) \\ &= \mathbf{\Phi}_0^T(x_1) \cdot \mathbf{q}_0 + \mathbf{\Phi}_1^T(x_1) \cdot \mathbf{q}_1 + \mathbf{\Phi}_s^T(x_1) \cdot \mathbf{q}_{fs}, \end{aligned} \tag{47}$$

where the vectors of undetermined constants $\mathbf{q}_0$ and $\mathbf{q}_{fs}$ are related to the source function $f(x_1)$ and its corresponding interpolation algebraical polynomial $f_s(x_1)$, and are determined by Eqs. (42) or (44), (46) and Eq. (17) respectively. However, the vector of undetermined constants $\mathbf{q}_1$ is to be determined by taking the equation above to satisfy the prescribed boundary conditions.

*2.6. Equivalent transformation of the solution*

For brevity, suppose that the displacement type boundary conditions are expressed by
$$\mathbf{q}_b = [\varphi(-a) \quad \varphi(a)]^{\mathrm{T}}. \tag{48}$$
Then substituting Eq. (47) into Eq. (48), we obtain
$$\mathbf{q}_1 = \mathbf{R}_f^{-1}\mathbf{q}_b - \mathbf{R}_f^{-1}\begin{bmatrix} \mathbf{\Phi}_0^{\mathrm{T}}(-a)\cdot\mathbf{q}_0 + \mathbf{\Phi}_s^{\mathrm{T}}(-a)\cdot\mathbf{q}_{fs} \\ \mathbf{\Phi}_0^{\mathrm{T}}(a)\cdot\mathbf{q}_0 + \mathbf{\Phi}_s^{\mathrm{T}}(a)\cdot\mathbf{q}_{fs} \end{bmatrix}, \tag{49}$$
where
$$\mathbf{R}_f = \begin{bmatrix} \mathbf{\Phi}_1^{\mathrm{T}}(-a) \\ \mathbf{\Phi}_1^{\mathrm{T}}(a) \end{bmatrix}. \tag{50}$$
Plugging back into Eq. (47) gives
$$\varphi(x_1) = \mathbf{\Phi}_1^{\mathrm{T}}(x_1)\cdot\mathbf{R}_f^{-1}\mathbf{q}_b + \mathbf{\Phi}_0^{\mathrm{T}}(x_1)\cdot\mathbf{q}_0 + \mathbf{\Phi}_s^{\mathrm{T}}(x_1)\cdot\mathbf{q}_{fs}$$
$$-\mathbf{\Phi}_1^{\mathrm{T}}(x_1)\cdot\mathbf{R}_f^{-1}\begin{bmatrix} \mathbf{\Phi}_0^{\mathrm{T}}(-a)\cdot\mathbf{q}_0 + \mathbf{\Phi}_s^{\mathrm{T}}(-a)\cdot\mathbf{q}_{fs} \\ \mathbf{\Phi}_0^{\mathrm{T}}(a)\cdot\mathbf{q}_0 + \mathbf{\Phi}_s^{\mathrm{T}}(a)\cdot\mathbf{q}_{fs} \end{bmatrix}. \tag{51}$$
We denote the vector of functions
$$\mathbf{\Phi}_b^{\mathrm{T}}(x_1) = \mathbf{\Phi}_1^{\mathrm{T}}(x_1)\cdot\mathbf{R}_f^{-1}, \tag{52}$$
and the functions
$$\varphi_b(x_1) = \mathbf{\Phi}_b^{\mathrm{T}}(x_1)\cdot\mathbf{q}_b, \tag{53}$$
$$\varphi_f(x_1) = \mathbf{\Phi}_0^{\mathrm{T}}(x_1)\cdot\mathbf{q}_0 + \mathbf{\Phi}_s^{\mathrm{T}}(x_1)\cdot\mathbf{q}_{fs} - \mathbf{\Phi}_b^{\mathrm{T}}(x_1)\cdot\begin{bmatrix} \mathbf{\Phi}_0^{\mathrm{T}}(-a)\cdot\mathbf{q}_0 + \mathbf{\Phi}_s^{\mathrm{T}}(-a)\cdot\mathbf{q}_{fs} \\ \mathbf{\Phi}_0^{\mathrm{T}}(a)\cdot\mathbf{q}_0 + \mathbf{\Phi}_s^{\mathrm{T}}(a)\cdot\mathbf{q}_{fs} \end{bmatrix}. \tag{54}$$
Then the Fourier series multiscale solution of the one-dimensional convection-diffusion-reaction equation can be rewritten as
$$\varphi(x_1) = \varphi_b(x_1) + \varphi_f(x_1). \tag{55}$$

# 3. One-dimensional numerical examples

As shown in Table 3, we have implemented the Fourier series multiscale method for the one-dimensional convection-diffusion-reaction equation with wide range of the computational parameters and boundary conditions, and presented an analysis of the effects both of the order of the interpolation algebraical polynomial of the source function and the derivation method for discrete equations on the convergence characteristics and approximation accuracy of the Fourier series multiscale solution. By this means, we optimize the settings of computational schemes of the Fourier series multiscale solution of the one-dimensional convection-diffusion-reaction equation. And with the optimized computational schemes, we investigate the multiscale characteristics of the one-dimensional convection-diffusion-reaction equation for varied computational parameters.

Table 3: Computational schemes of the one-dimensional convection-diffusion-reaction equation.

| Configuration of problem | | Computational scheme | |
| --- | --- | --- | --- |
| Boundary condition | Computational parameter $(P_e, D_a)$ | Order of interpolation algebraical polynomial of source function | Derivation method for discrete equations |
| DD | (3, 90) | $N_1^s = 0$ | FCCM |
| DN | (1, 30) | $N_1^s = 1$ | CM |
|  | (30, 1) | $N_1^s = 2$ |  |
|  | (200, -1) |  |  |

*3.1. Convergence characteristics*

As shown in Table 4, we perform four comparative convergence experiments in series for detailed analysis of the influences of some possible factors on convergence characteristics of the Fourier series multiscale solution. These possible factors are extended to involve in the one-dimensional convection-diffusion-reaction equation (computational parameters and boundary conditions) and its Fourier series multiscale solution (the order of the interpolation algebraical polynomial of the source function and the derivation method for discrete equations). We adopt the computational scheme specified in the first column of Table 3 as reference, and then adjustment of the single factors (see the corresponding rows in Table 3) leads to the four comparative convergence experiments. For instance, in the first comparative convergence experiment, we start out with the reference computational scheme, and change the derivation method for discrete equations from the Fourier coefficient comparison method to the collocation method; in the second comparative convergence experiment, we start out with the reference computational scheme, and change the order of the interpolation algebraical polynomial of the source function from zero to one and two successively; in the third comparative convergence experiment, we start out with the reference computational scheme, and change computational parameters of the one-dimensional convection-diffusion-reaction equation from strong reaction type to reaction-dominated type, convection-dominated type and strong convection type successively; and in the fourth comparative convergence experiment, we start out with the reference computational scheme, and change the boundary condition of the one-dimensional convection-diffusion-reaction equation from the DD boundary condition to the DN boundary condition.

In the comparative convergence experiments, the source function remains unchanged, with the following form of the third order algebraical polynomial

$$f(x_1) = 10^3 + 2 \times 10^3 \frac{x_1}{a} + 5 \times 10^3 (\frac{x_1}{a})^2 + 10^4 (\frac{x_1}{a})^3, \quad x_1 \in [-a, a]. \tag{56}$$

Therefore, when we adjust the order of the interpolation algebraical polynomial of the source function to three, it follows that the vector of Fourier coefficients

$$\mathbf{q}_0 = \mathbf{0}, \tag{57}$$

and accordingly the Fourier series multiscale solution given by Eq. (47) yields the exact solution of the one-dimensional convection-diffusion-reaction equation.

Table 4: Comparative convergence experiments for one-dimensional convection-diffusion-reaction equation.

| Comparative convergence experiment | No. | Boundary condition | Computational parameter $(P_e, D_a)$ | Order of interpolation algebraical polynomial of source function | Derivation method for discrete equations |
|---|---|---|---|---|---|
| 1 | a | DD $\varphi(-a)=1$, $\varphi(a)=0$ | (3, 90) | $N_1^s = 0$ | FCCM |
|   | b |   |   |   | CM |
| 2 | a | DD $\varphi(-a)=1$, $\varphi(a)=0$ | (3, 90) | $N_1^s = 0$ | FCCM |
|   | b |   |   | $N_1^s = 1$ |   |
|   | c |   |   | $N_1^s = 2$ |   |
| 3 | a | DD $\varphi(-a)=1$, $\varphi(a)=0$ | (3, 90) | $N_1^s = 0$ | FCCM |
|   | b |   | (1, 30) |   |   |
|   | c |   | (30, 1) |   |   |
|   | d |   | (200, -1) |   |   |
| 4 | a | DD $\varphi(-a)=1$, $\varphi(a)=0$ | (3, 90) | $N_1^s = 0$ | FCCM |
|   | b | DN $\varphi(-a)=1$, $\varphi^{(1)}(a)=1$ |   |   |   |

1. General convergence characteristics

As to the reference computational scheme given in Table 3, we truncate the composite Fourier series of the Fourier series multiscale solution successively with the first 2, 3, 5, 10, 20, 30 and 40 terms and then we compute the internal approximation errors and boundary approximation errors of the function $\varphi(x_1)$ and its first and second order derivatives. Some of the results are shown in Figure 1.

With the increase of the number of truncated terms of the composite Fourier series, the evolution of computed values of the indexes of approximation errors is exhibited as the following:

a. The internal approximation errors of $\varphi(x_1)$ decrease rapidly as the number of truncated terms of the composite Fourier series increases. When the number of truncated terms equals to 10, the approximation error is already about 1.0E-3 and shows a trend of sustained and rapid decrease. And when the number of truncated terms equals to 40, the approximation error decreases further by 2 orders of magnitude.

b. The boundary approximation errors of $\varphi(x_1)$ remain extremely small as the number of truncated terms of the composite Fourier series increases. The approximation errors are all below 1.0E-16.

c. The internal approximation errors of $\varphi^{(1)}(x_1)$ decrease rapidly as the number of

truncated terms of the composite Fourier series increases. When the number of truncated terms equals to 10, the approximation error is already about 1.0E-2 and shows a trend of sustained decrease. And when the number of truncated terms equals to 40, the approximation error decreases further by 1 order of magnitude.

d. The boundary approximation errors of $\varphi^{(1)}(x_1)$ decrease slowly as the number of truncated terms of the composite Fourier series increases. When the number of truncated terms equals to 10, the approximation error is already about 1.0E-1 and shows a trend of sustained decrease. And when the number of truncated terms equals to 40, the approximation error decreases further by 1 order of magnitude.

e. The internal approximation errors of $\varphi^{(2)}(x_1)$ decrease rapidly as the number of truncated terms of the composite Fourier series increases. When the number of truncated terms equals to 10, the approximation error is already about 1.0E-2 and shows a trend of slow-down decrease. And when the number of truncated terms equals to 40, the approximation error decreases further by 0.5 order of magnitude.

f. The boundary approximation errors of $\varphi^{(2)}(x_1)$ remain almost unchanged as the number of truncated terms of the composite Fourier series increases. The approximation errors are between 0.5 and 1.0.

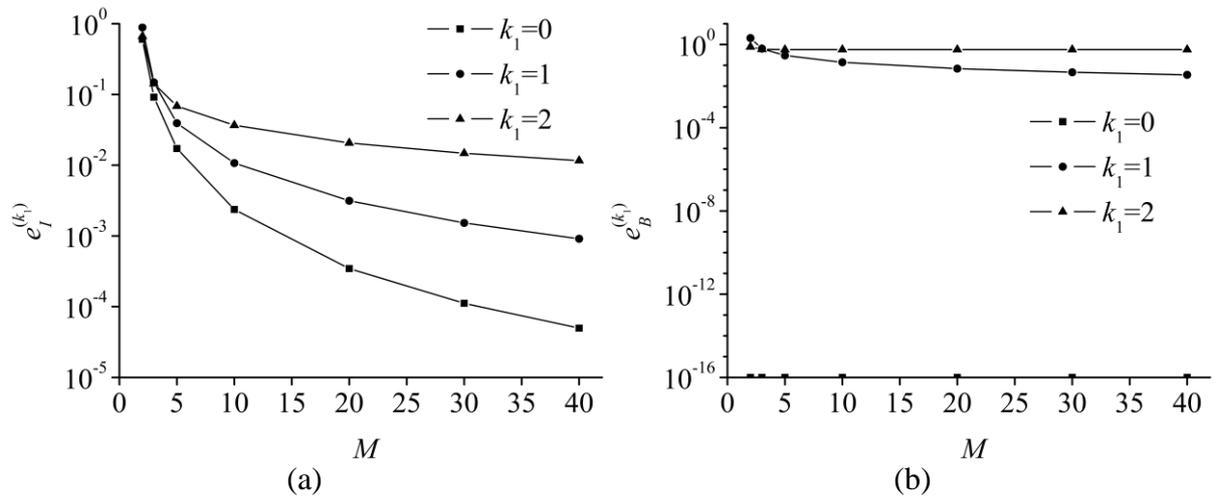

Figure 1: The convergence characteristics of Fourier series multiscale solutions of one-dimensional convection-diffusion-reaction equation:
(a) $e_I^{(k_1)}(\varphi_M)$-$M$ curves, (b) $e_B^{(k_1)}(\varphi_M)$-$M$ curves.

2. Comparative convergence experiment 1

When the derivation method for discrete equations changes from the Fourier coefficient comparison method (FCCM) to the collocation method (CM), the evolution of convergence characteristics of the Fourier series multiscale solution is illustrated in Figure 2. We make a brief analysis as the following:

a. The convergence of the composite Fourier series of $\varphi(x_1)$, $\varphi^{(1)}(x_1)$ and $\varphi^{(2)}(x_1)$ remain almost unchanged within the solution domain, namely the composite Fourier series all converge well within the solution domain. However, when the number of truncated terms of the composite Fourier series equals to 40, the corresponding approximation errors increase by 1 to 2 orders of magnitude.

b. The convergence of the composite Fourier series of $\varphi(x_1)$, $\varphi^{(1)}(x_1)$ and $\varphi^{(2)}(x_1)$ remain almost unchanged on the boundary of the solution domain, namely the composite Fourier series of $\varphi(x_1)$, $\varphi^{(1)}(x_1)$ converge well on the boundary, and the composite Fourier series of $\varphi^{(2)}(x_1)$ do not converge as well on the boundary.

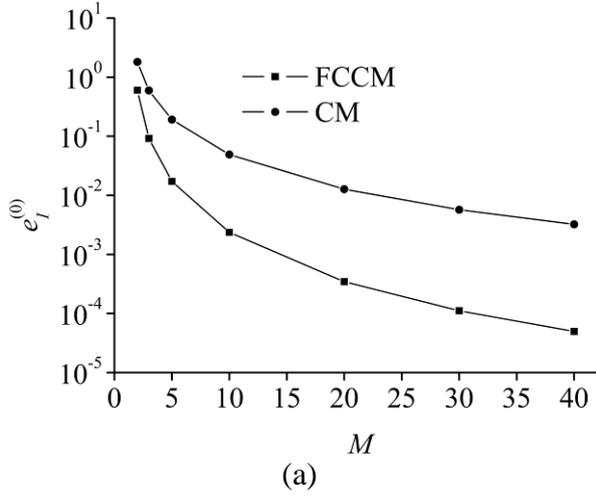
(a)

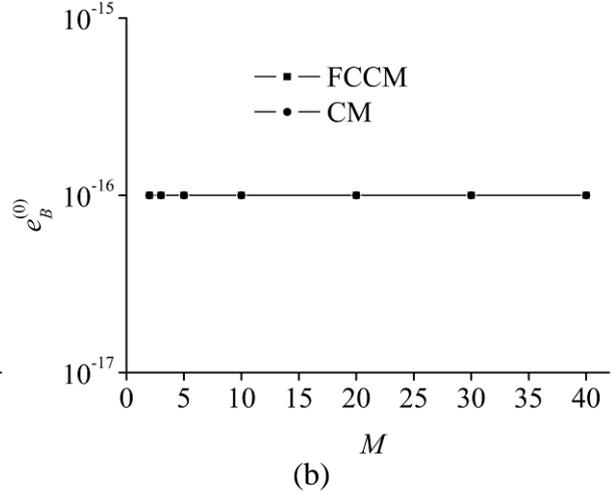
(b)

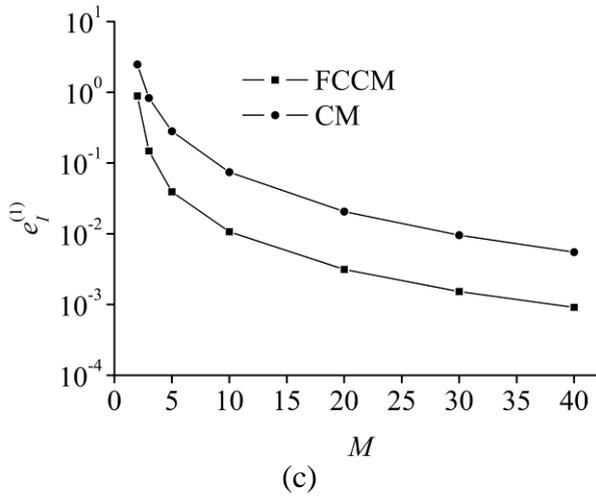
(c)

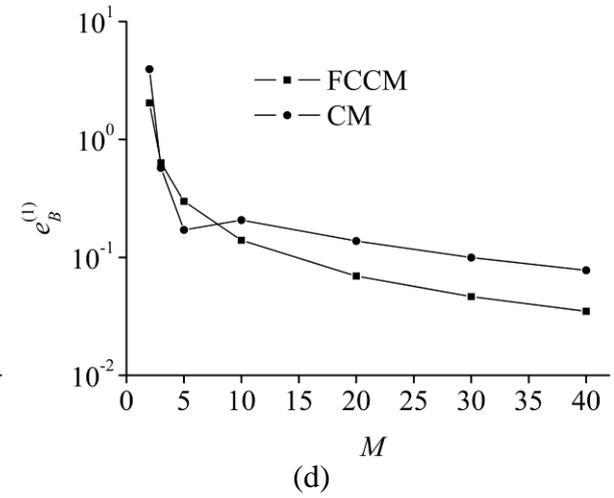
(d)

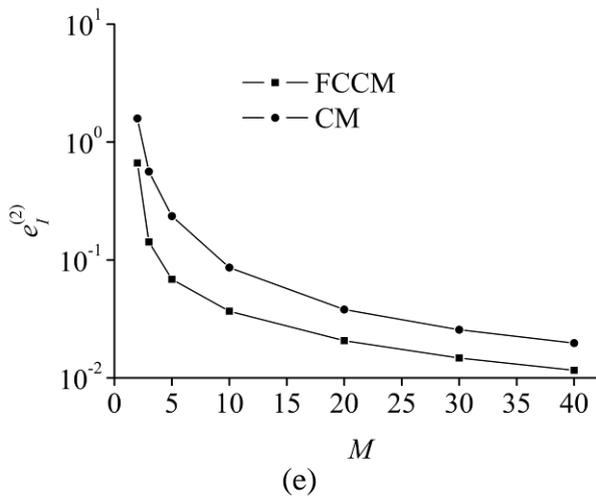
(e)

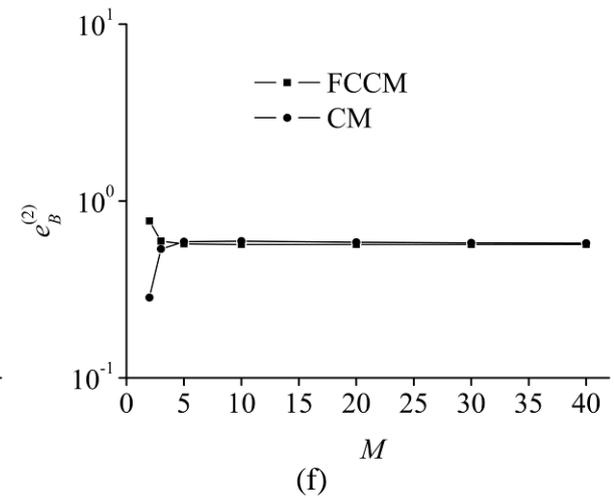
(f)

Figure 2: Convergence comparison of the Fourier series multiscale solutions with different derivation methods for discrete equations (FCCM and CM):

(a) $e_I^{(0)}(\varphi_M)$-$M$ curves, (b) $e_B^{(0)}(\varphi_M)$-$M$ curves, (c) $e_I^{(1)}(\varphi_M)$-$M$ curves,

(d) $e_B^{(1)}(\varphi_M)$-$M$ curves, (e) $e_I^{(2)}(\varphi_M)$-$M$ curves, (f) $e_B^{(2)}(\varphi_M)$-$M$ curves.

3. Comparative convergence experiment 2

When the order of the interpolation algebraical polynomial of the source function changes from zero to one and two successively, the evolution of convergence characteristics of the Fourier series multiscale solution is illustrated in Figure 3. We make a brief analysis as the following:

a. The convergence of the composite Fourier series of $\varphi(x_1)$, $\varphi^{(1)}(x_1)$ and $\varphi^{(2)}(x_1)$ remain almost unchanged within the solution domain, namely the composite Fourier series all converge well within the solution domain. However, when the number of truncated terms of the composite Fourier series equals to 40, the corresponding approximation errors decrease respectively by 1 to 3, 1 to 3, and 3 to 4 orders of magnitude.

b. The convergence of the composite Fourier series of $\varphi(x_1)$ and $\varphi^{(1)}(x_1)$ remain almost unchanged on the boundary of the solution domain, namely the composite Fourier series all converge well on the boundary. Moreover, the approximation errors of $\varphi(x_1)$ remain extremely small (all below 1.0E-16) as the number of truncated terms of the composite Fourier series increases. And, when the number of truncated terms of the composite Fourier series equals to 40, the approximation errors of $\varphi^{(1)}(x_1)$ decrease respectively by 3 and 4 orders of magnitude.

c. The convergence of the composite Fourier series of $\varphi^{(2)}(x_1)$ is improved significantly on the boundary of the solution domain. When the number of truncated terms of the composite Fourier series equals to 40, the corresponding approximation errors decrease respectively by 2 and 6 orders of magnitude.

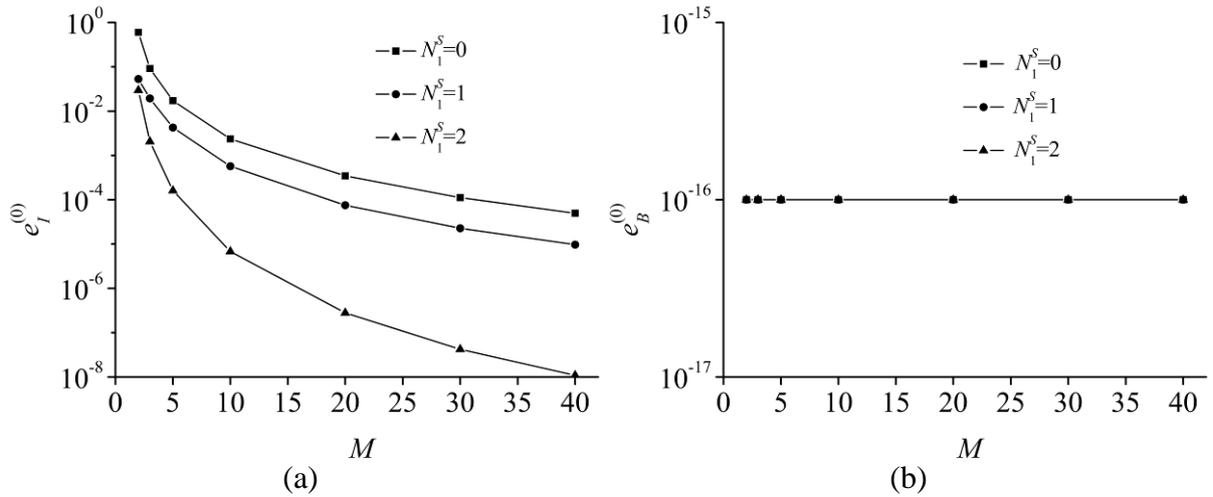

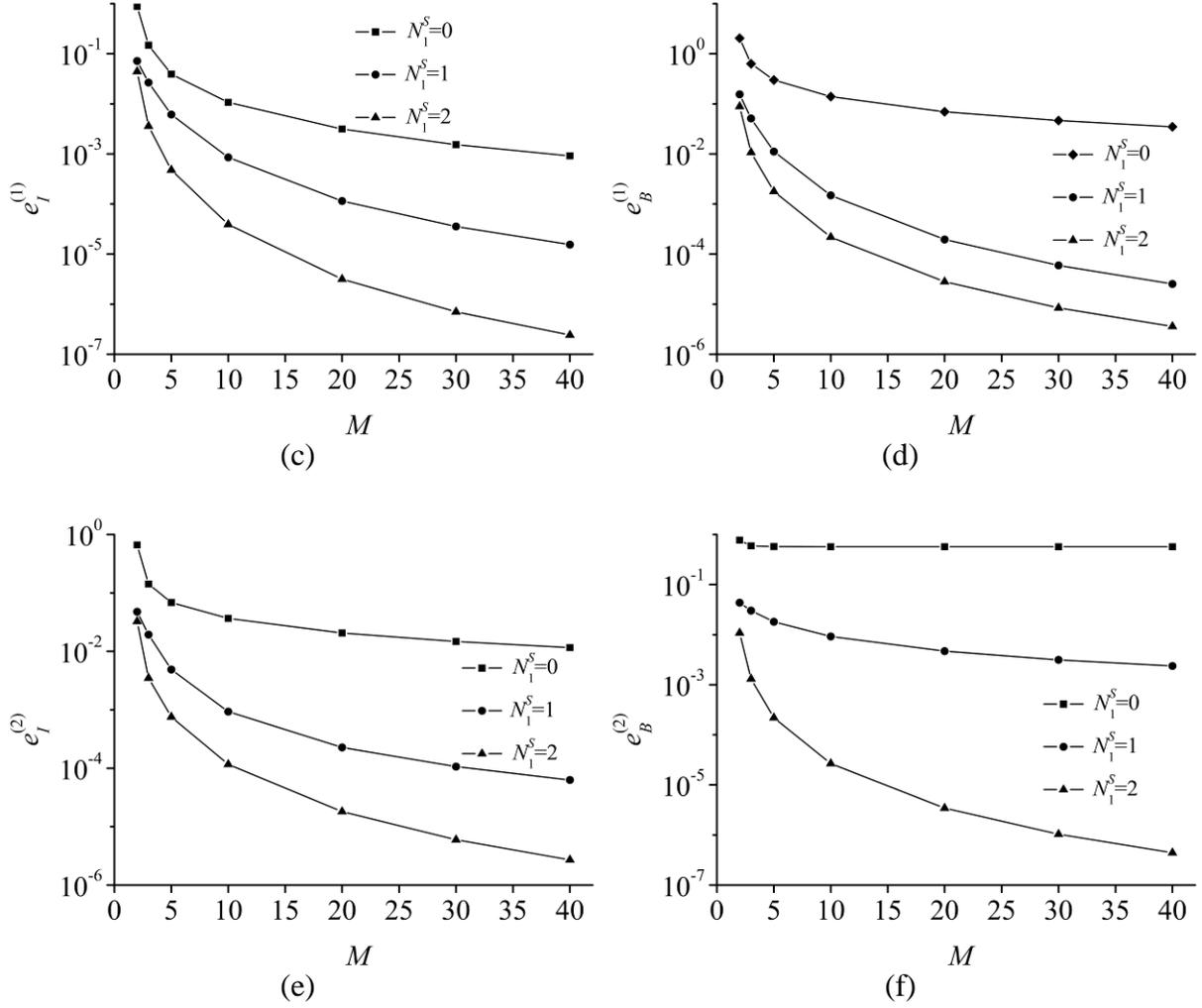

Figure 3: Convergence comparison of the Fourier series multiscale solutions with different orders of interpolation polynomials of source function:
(a) $e_I^{(0)}(\varphi_M)$-$M$ curves, (b) $e_B^{(0)}(\varphi_M)$-$M$ curves, (c) $e_I^{(1)}(\varphi_M)$-$M$ curves,
(d) $e_B^{(1)}(\varphi_M)$-$M$ curves, (e) $e_I^{(2)}(\varphi_M)$-$M$ curves, (f) $e_B^{(2)}(\varphi_M)$-$M$ curves.

4. Comparative convergence experiment 3

When the computational parameters of the one-dimensional convection-diffusion-reaction equation change from strong reaction type ($P_e = 3$, $D_a = 90$) to reaction-dominated type ($P_e = 1$, $D_a = 30$), convection-dominated type ($P_e = 30$, $D_a = 1$) and strong convection type ($P_e = 200$, $D_a = -1$) successively, the evolution of convergence characteristics of the Fourier series multiscale solution is illustrated in Figure 4. We make a brief analysis as the following:

a. The convergence of the composite Fourier series of $\varphi(x_1)$, $\varphi^{(1)}(x_1)$ and $\varphi^{(2)}(x_1)$ remain almost unchanged within the solution domain, namely the composite Fourier series all converge well within the solution domain. However, when the number of truncated terms of the composite Fourier series equals to 40, the corresponding approximation errors are of difference in quantity up to 3 orders of magnitude.

b. The convergence of the composite Fourier series of $\varphi(x_1)$ and $\varphi^{(1)}(x_1)$ remain

almost unchanged on the boundary of the solution domain, namely the composite Fourier series all converge well on the boundary. Moreover, the approximation errors of $\varphi(x_1)$ remain extremely small (all below 1.0E-13) as the number of truncated terms of the composite Fourier series increases. And, when the number of truncated terms of the composite Fourier series equals to 40, the corresponding approximation errors of $\varphi^{(1)}(x_1)$ are of difference in quantity up to 1 order of magnitude.

c. The convergence of the composite Fourier series of $\varphi^{(2)}(x_1)$ are of marked difference on the boundary of the solution domain. Of the four sets of computational parameters, the composite Fourier series converge well for the case of $P_e = 200$ and $D_a = -1$, and do not converge as well for other cases.

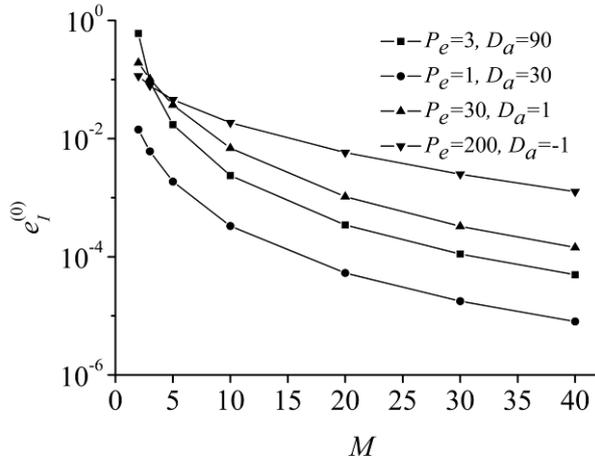

(a)

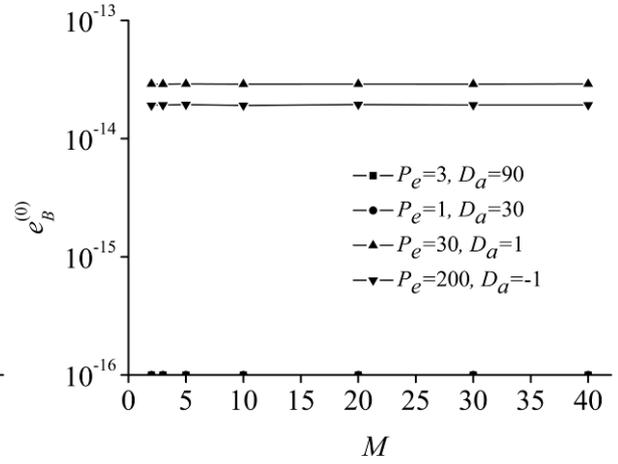

(b)

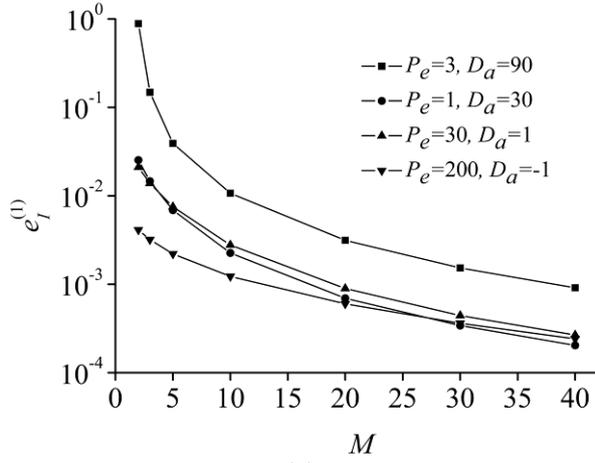

(c)

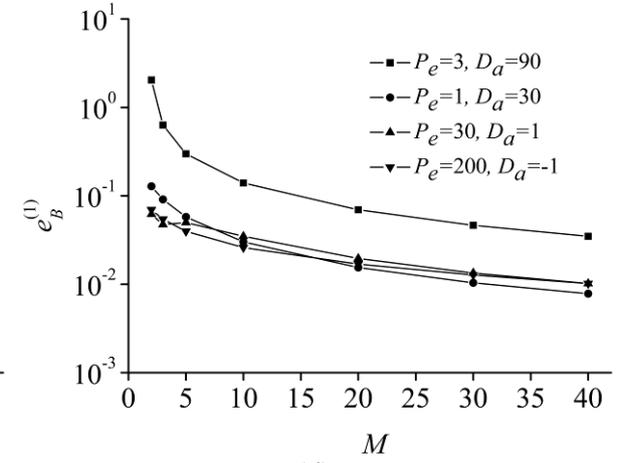

(d)

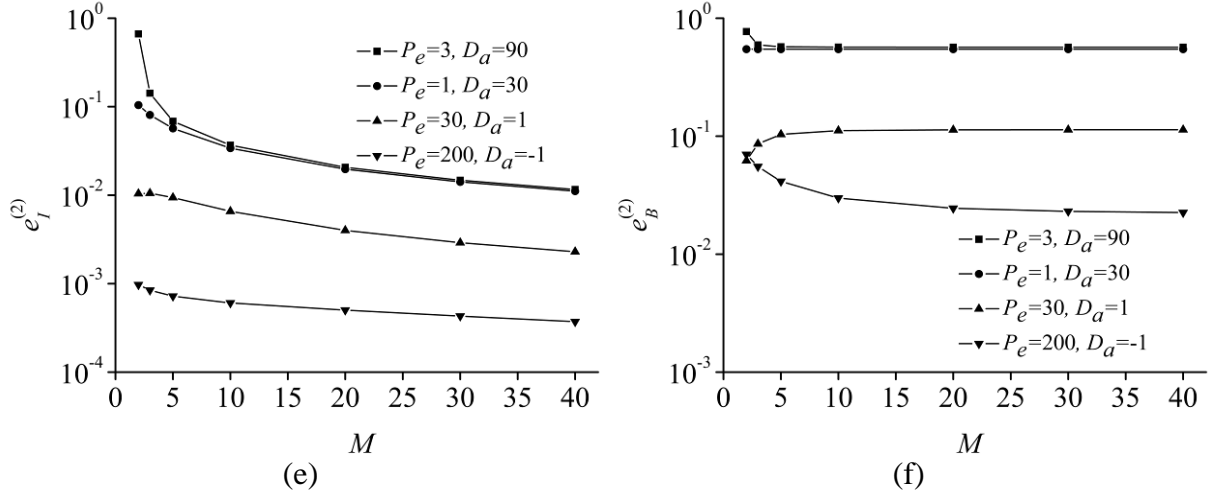

Figure 4: Convergence comparison of the Fourier series multiscale solutions with different computational parameters $P_e$ and $D_a$:

(a) $e_I^{(0)}(\varphi_M)$-$M$ curves, (b) $e_B^{(0)}(\varphi_M)$-$M$ curves, (c) $e_I^{(1)}(\varphi_M)$-$M$ curves,

(d) $e_B^{(1)}(\varphi_M)$-$M$ curves, (e) $e_I^{(2)}(\varphi_M)$-$M$ curves, (f) $e_B^{(2)}(\varphi_M)$-$M$ curves.

5. Comparative convergence experiment 4

When the boundary condition of the one-dimensional convection-diffusion-reaction equation changes from the DD boundary condition to the DN boundary condition, the evolution of convergence characteristics of the Fourier series multiscale solution is illustrated in Figure 5. We make a brief analysis as the following:

a. The convergence of the composite Fourier series of $\varphi(x_1)$, $\varphi^{(1)}(x_1)$ and $\varphi^{(2)}(x_1)$ remain almost unchanged within the solution domain, namely the composite Fourier series all converge well within the solution domain. However, when the number of truncated terms of the composite Fourier series equals to 40, the corresponding approximation errors respectively increase by 2 orders of magnitude, increase by 1 order of magnitude or remain almost unchanged.

b. The convergence of the composite Fourier series of $\varphi(x_1)$ is worsened slightly on the boundary of the solution domain. When the number of truncated terms of the composite Fourier series equals to 40, the corresponding approximation errors increase from 1.0E-16 to 1.0E-2.

c. The convergence of the composite Fourier series of $\varphi^{(1)}(x_1)$ is improved significantly on the boundary of the solution domain. When the number of truncated terms of the composite Fourier series equals to 40, the corresponding approximation errors decrease by 1 order of magnitude.

d. The convergence of the composite Fourier series of $\varphi^{(2)}(x_1)$ remains almost unchanged on the boundary of the solution domain, namely the composite Fourier series of $\varphi^{(2)}(x_1)$ do not converge as well on the boundary.

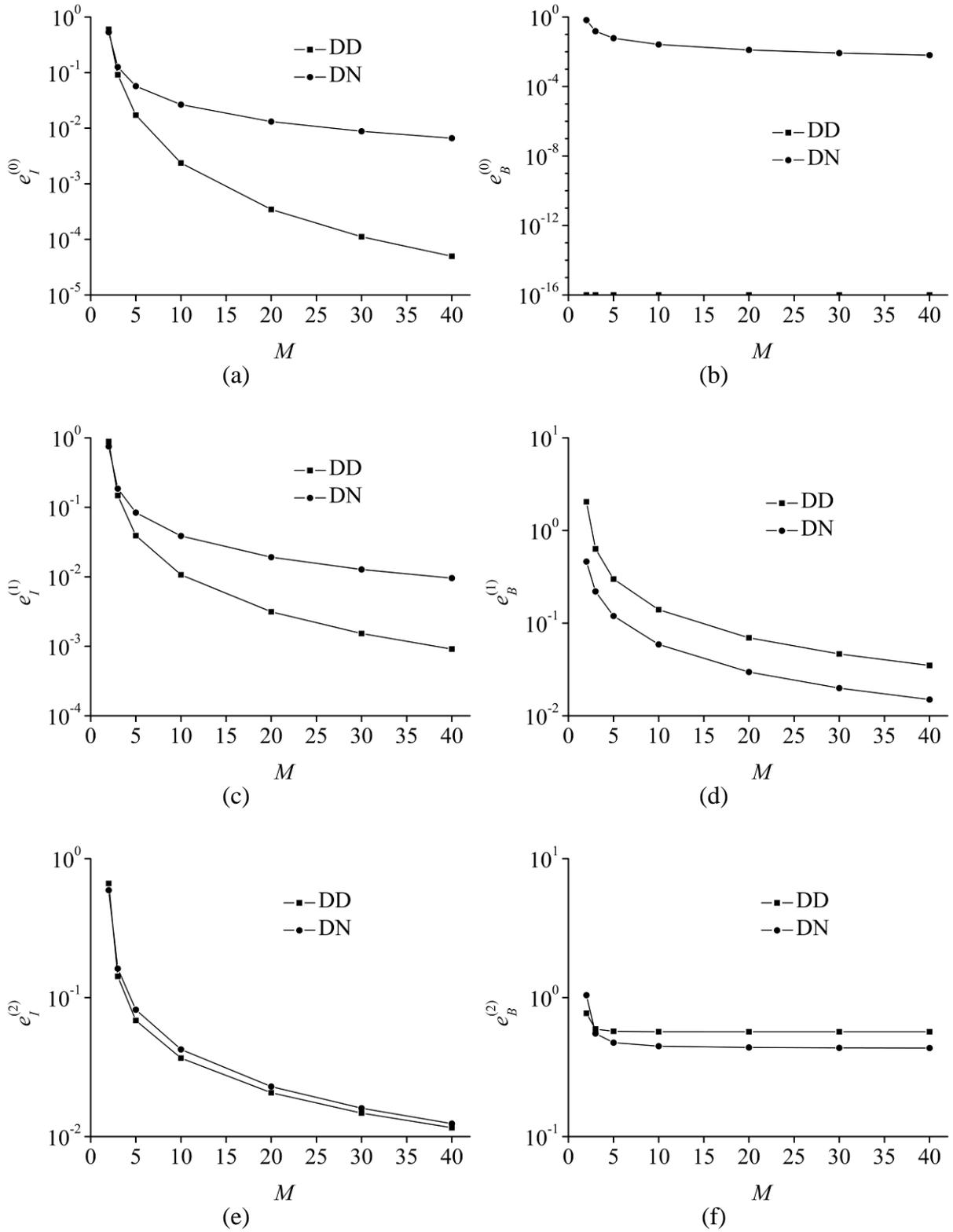

Figure 5: Convergence comparison of the Fourier series multiscale solutions with different boundary conditions:
(a) $e_I^{(0)}(\varphi_M)$-$M$ curves, (b) $e_B^{(0)}(\varphi_M)$-$M$ curves, (c) $e_I^{(1)}(\varphi_M)$-$M$ curves,
(d) $e_B^{(1)}(\varphi_M)$-$M$ curves, (e) $e_I^{(2)}(\varphi_M)$-$M$ curves, (f) $e_B^{(2)}(\varphi_M)$-$M$ curves.

6. Brief summary

Therefore, the convergence characteristics of the Fourier series multiscale solution of the one-dimensional convection-diffusion-reaction equation are summarized as follows:

a. The Fourier series multiscale solution has good convergence in general. And especially, the composite Fourier series of $\varphi(x_1)$ and its first order and second order derivatives converge well within the solution domain, the composite Fourier series of $\varphi(x_1)$ and its first order derivative converge well on the boundary of the solution domain, and the composite Fourier series of the second order derivative of $\varphi(x_1)$ does not converge as well on the boundary of the solution domain.

b. The influence factors, as investigated for the convergence characteristics of the Fourier series multiscale solution, display varied behaviors. For instance, the adjustment of the derivation method for discrete equations brings about few effects on the convergence of the Fourier series multiscale solution, the increase of the order of the interpolation algebraical polynomial of the source function significantly improves the convergence of the Fourier series multiscale solution, and the adjustment of the computational parameters and boundary conditions has marked influences on the convergence characteristics of the Fourier series multiscale solution.

*3.2. Computational efficiency*

To verify the computational efficiency of the Fourier series multiscale method, we compare it with the bubble function method (see [18]) for the one-dimensional convection-diffusion-reaction equation. We carry out this comparison by means of the change of the computational parameters from the strong reaction type successively to the reaction-dominated type, the convection-dominated type and the strong convection type, as specified in the third comparative convergence experiment of Table 4. And in this process, we employ the source functions in the form of the constant

$$f(x_1) = 10^3, \quad x_1 \in [-a, a], \tag{58}$$

or the first order algebraical polynomial

$$f(x_1) = 10^3 + 2 \times 10^3 \frac{x_1}{a}, \quad x_1 \in [-a, a]. \tag{59}$$

In the specific comparative experiment of computational efficiency, we truncate the composite Fourier series with the first 40 terms for the Fourier series multiscale method, and partition the solution interval into 100 uniform meshes for the bubble function method. And thus, for the eight different cases that are taken into consideration, we list in Table 5 some performance indexes, such as the numbers of undetermined constants, the computing time (codes are running on the same computer), and the overall computational errors (the definition refers to Section 5.1 in [40], and the number of sampling points is 10001), of the Fourier series multiscale method and the bubble function method. Herein we make a brief analysis of the computational efficiency of these two different types of multiscale methods:

1. The number of undetermined constants involved in the Fourier series multiscale method is 83, and the number of undetermined constants involved in the bubble function method is 101. Therefore, the computing scale of the Fourier series multiscale method is slightly smaller than that of the bubble function method.

2. The code of the Fourier series multiscale method runs for 2 seconds on the computer, and the code of the bubble function method runs for 49 seconds on the computer. Therefore, the computing time of the Fourier series multiscale method is remarkably shorter than that of the bubble function method.

3. The Fourier series multiscale method seems to be insensitive to the change of computational parameters and form of the source function, and is applicable to all the eight different cases in the comparative experiment of computational efficiency. However, the bubble function method is markedly influenced by the change of computational parameters and form of the source function, and is only applicable to the six different cases in the comparative experiment of computational efficiency.

4. The computational errors of the Fourier series multiscale method are between 1.0E-15 and 1.0E-2, and the computational errors of the bubble function method are between 1.0E-4 and 7.0E-1. Therefore, the computational accuracy of the Fourier series multiscale method is superior to that of the bubble function method.

Hence, it is verified that the Fourier series multiscale method excels the bubble function method in computational efficiency and has shown considerable potential to be a practical, feasible, efficient and highly accurate method.

Table 5: Comparison between computational efficiency of Fourier series multiscale method and bubble function method

| Method | Number of constants | Computing time | Source function | Computational Parameter $(P_e, D_a)$ | Computational error $\varphi$ | $\varphi^{(1)}$ | $\varphi^{(2)}$ |
|---|---|---|---|---|---|---|---|
| Fourier series multiscale method | 83 | 2s | Constant | (3, 90) | 1.0849E-13 | 9.8395E-14 | 1.1072E-13 |
| | | | | (1, 30) | 1.8423E-13 | 1.9183E-13 | 1.9945E-13 |
| | | | | (30, 1) | 1.6526E-13 | 2.0002E-14 | 1.2077E-14 |
| | | | | (200, -1) | 1.8946E-13 | 3.1784E-15 | 1.6337E-15 |
| Fourier series multiscale method | 83 | 2s | First order algebraical polynomial | (3, 90) | 3.1502E-05 | 8.6717E-04 | 1.2883E-02 |
| | | | | (1, 30) | 3.4250E-06 | 1.3947E-04 | 7.7467E-03 |
| | | | | (30, 1) | 4.1961E-05 | 9.9038E-05 | 9.5799E-04 |
| | | | | (200, -1) | 5.6963E-04 | 1.1687E-04 | 1.8354E-04 |
| Bubble function method | 101 | 49s | Constant | (3, 90) | 2.2078E-04 | 3.8535E-03 | 1.0692E-01 |
| | | | | (1, 30) | 1.0400E-04 | 8.3772E-03 | 7.0050E-01 |
| | | | | (30, 1) | 1.6512E-04 | 1.4376E-03 | 1.8817E-02 |
| | | | | (200, -1) | 2.2520E-03 | 3.3312E-03 | 6.9497E-03 |
| Bubble function method | 101 | 49s | First order algebraical polynomial | (3, 90) | 2.9823E-02 | 5.9332E-02 | 4.0030E-01 |
| | | | | (1, 30) | 3.9407E-02 | 1.2020E-01 | 2.6545E-01 |
| | | | | (30, 1) | Failure | | |
| | | | | (200, -1) | Failure | | |

## 3.3. Multiscale characteristics

As shown in Figure 6(a), the source function is of the form
$$f(x_1) = 1000\delta(x_1 - 0.5), \quad x_1 \in [0,1]. \tag{60}$$
We take this quasi-Green's function problem as an example and investigate the multiscale characteristics of the one-dimensional convection-diffusion-reaction equation with varying computational parameters.

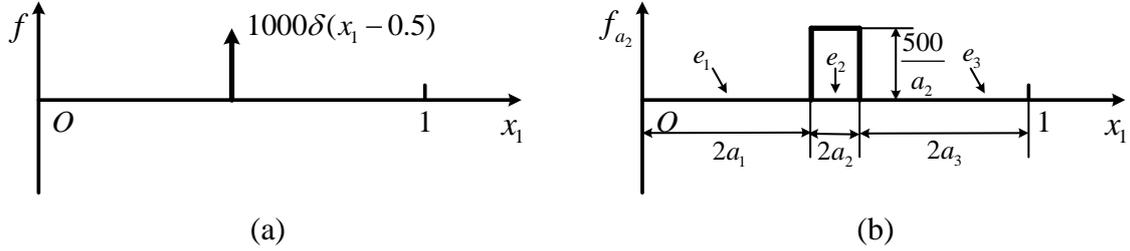

Figure 6: Descriptions of source function: (a) theoretical description, (b) approximate description.

For the above quasi-Green's function problem, we present in Table 6 two computational schemes such as the whole interval solution and the subinterval solution.

1. Translating the one-dimensional solution interval from $[0,1]$ to $[-1/2,1/2]$ and setting $a = 1/2$, we solve the convection-diffusion-reaction equation directly by the Fourier series multiscale method without the introduction of interpolation algebraical polynomial of the source function.

2. As shown in Figure 6(b), the source function is actually the limit of the sequence of functions $\{f_{a_2}(x_1)\}_{0<a_2<0.5}$ as the parameter $a_2$ approaches 0 from above. We solve the convection-diffusion-reaction equation respectively on the subintervals $e_1 = [0, 2a_1]$, $e_2 = [2a_1, 2a_1 + 2a_2]$ and $e_3 = [2a_1 + 2a_2, 1]$ by the Fourier series multiscale method with the introduction of first order interpolation algebraical polynomial of the source function. It is obvious that the Fourier series multiscale method yields accurate solution on each subinterval.

Table 6: Computational scheme for quasi-Green's function problem in one-dimensional convection-diffusion-reaction equation.

| Description of the problem | | Solution schemes | |
| --- | --- | --- | --- |
| Boundary conditions | Computational parameter $(P_e, D_a)$ | Whole interval solution | Subinterval solution |
| DD $\varphi(0) = 1, \ \varphi(1) = 0$ | (3, 90) | $N_1^s = 0$ FCCM | $N_1^s = 1$ $\mathbf{q}_0 = \mathbf{0}$ |
| | (1, 30) | | |
| | (30, 1) | | |
| | (200, -1) | | |

In the computational scheme of subinterval solution, we denote the endpoints of the subintervals by 1, 2, 3 and 4 respectively, and accordingly we denote the values of the solution function at these endpoints by $\varphi_1$, $\varphi_2$, $\varphi_3$ and $\varphi_4$ respectively. Using the displacement type boundary conditions at the first and fourth endpoints and the subinterval connection conditions at the second and third endpoints, we obtain

$$\left.\begin{array}{c}\varphi_1 = 1 \\ \Phi_{b,e_1}^{(1)\mathrm{T}}(a_1)[\varphi_1\ \varphi_2]^{\mathrm{T}} + \varphi_{f,e_1}^{(1)}(a_1) = \Phi_{b,e_2}^{(1)\mathrm{T}}(-a_2)[\varphi_2\ \varphi_3]^{\mathrm{T}} + \varphi_{f,e_2}^{(1)}(-a_2) \\ \Phi_{b,e_2}^{(1)\mathrm{T}}(a_2)[\varphi_2\ \varphi_3]^{\mathrm{T}} + \varphi_{f,e_2}^{(1)}(a_2) = \Phi_{b,e_3}^{(1)\mathrm{T}}(-a_3)[\varphi_3\ \varphi_4]^{\mathrm{T}} + \varphi_{f,e_3}^{(1)}(-a_3) \\ \varphi_4 = 0\end{array}\right\}. \qquad (61)$$

Hence, we determine $\varphi_2$ and $\varphi_3$, the values of the solution function at the specific endpoints. And further, by using Eq. (55), we determine the solution function on each subinterval and thus the solution function on the solution interval $[0,1]$.

Of the two computational schemes, the whole interval solution and the subinterval solution correspond respectively to the convergence process where $M$, the number of truncated terms of the composite Fourier series, approaches $\infty$, and the convergence process where $a_2$, the size parameter of the subinterval $e_2$, approaches 0 from above. As shown in Figures 7-10, different computational parameters lead to clear differences in the modes of the convergence process. And for some specific computational parameters, the solution of the one-dimensional convection-diffusion-reaction equation exhibits multiscale phenomena, such as sharp gradients and singularities, in the solution interval.

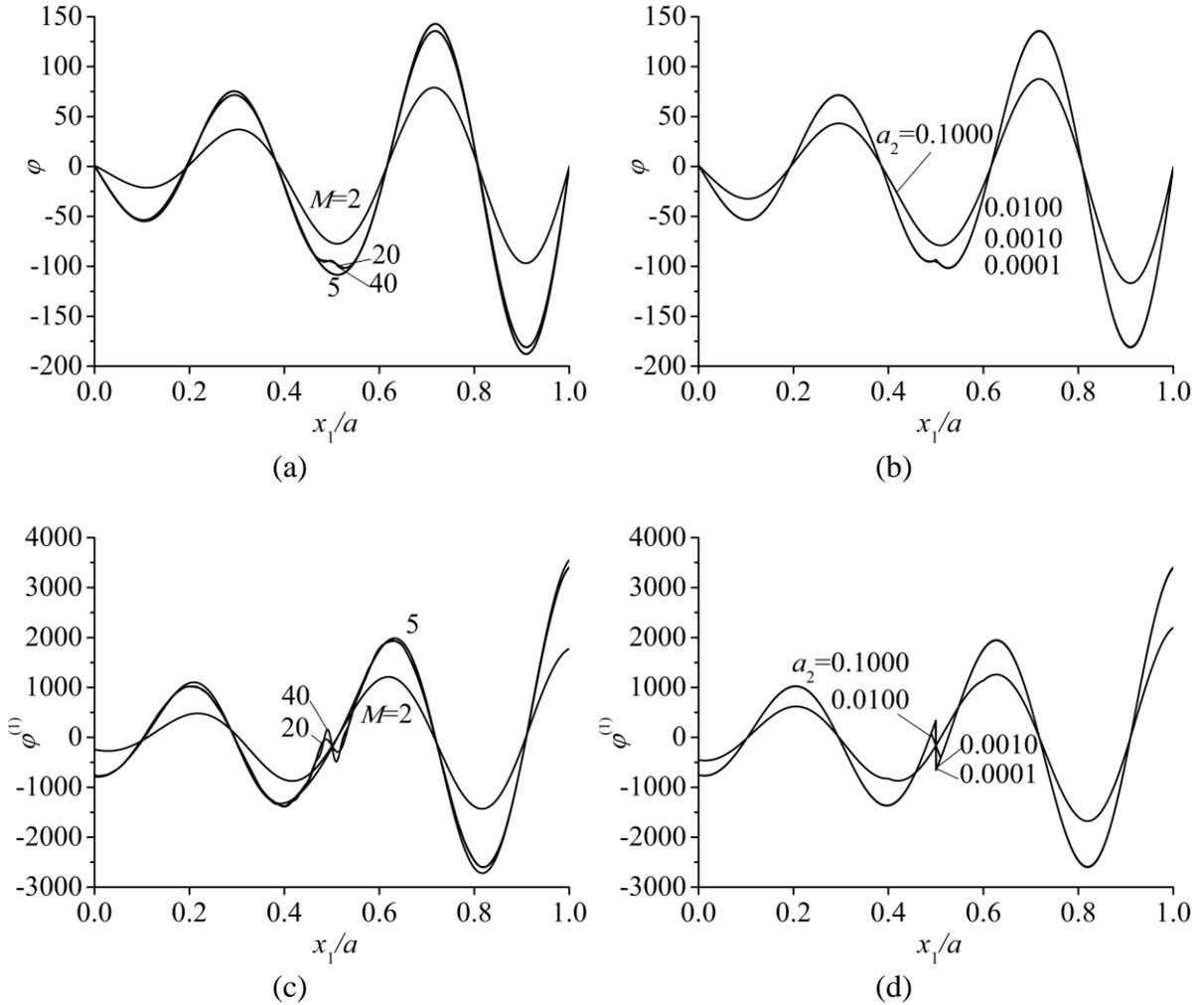

(a)  (b)  (c)  (d)

Figure 7: Two convergence processes for the whole interval solution and the subinterval solution with computational parameters $P_e = 3$ and $D_a = 90$:

(a) $\varphi$ curves for the whole interval solution with different terms of $M$,

(b) $\varphi$ curves for the subinterval solution with different values of $a_2$,

(c) $\varphi^{(1)}$ curves for the whole interval solution with different terms of $M$,

(d) $\varphi^{(1)}$ curves for the subinterval solution with different values of $a_2$.

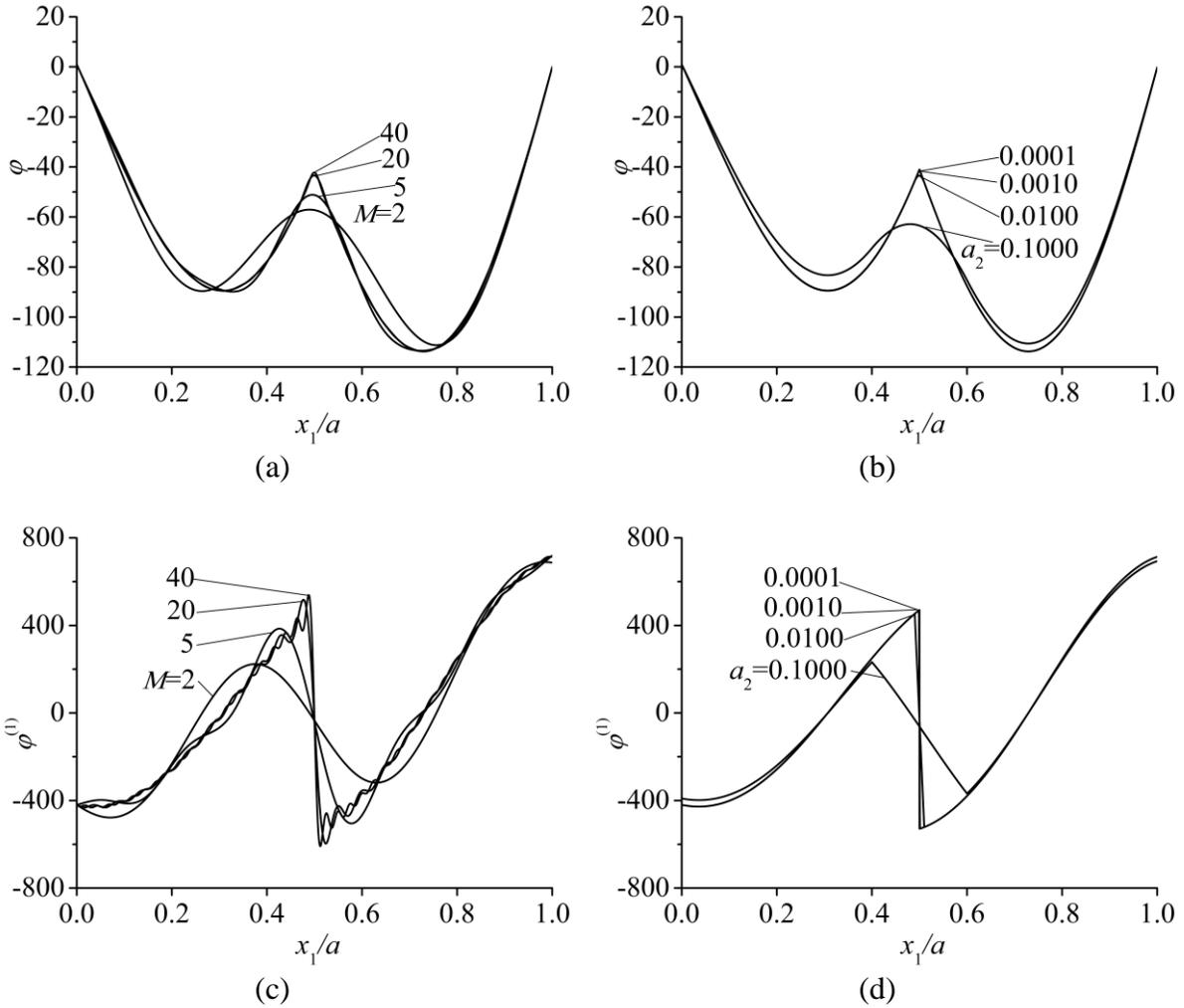

Figure 8: Two convergence processes for the whole interval solution and the subinterval solution with computational parameters $P_e = 1$ and $D_a = 30$:

(a) $\varphi$ curves for the whole interval solution with different terms of $M$,

(b) $\varphi$ curves for the subinterval solution with different values of $a_2$,

(c) $\varphi^{(1)}$ curves for the whole interval solution with different terms of $M$,

(d) $\varphi^{(1)}$ curves for the subinterval solution with different values of $a_2$.

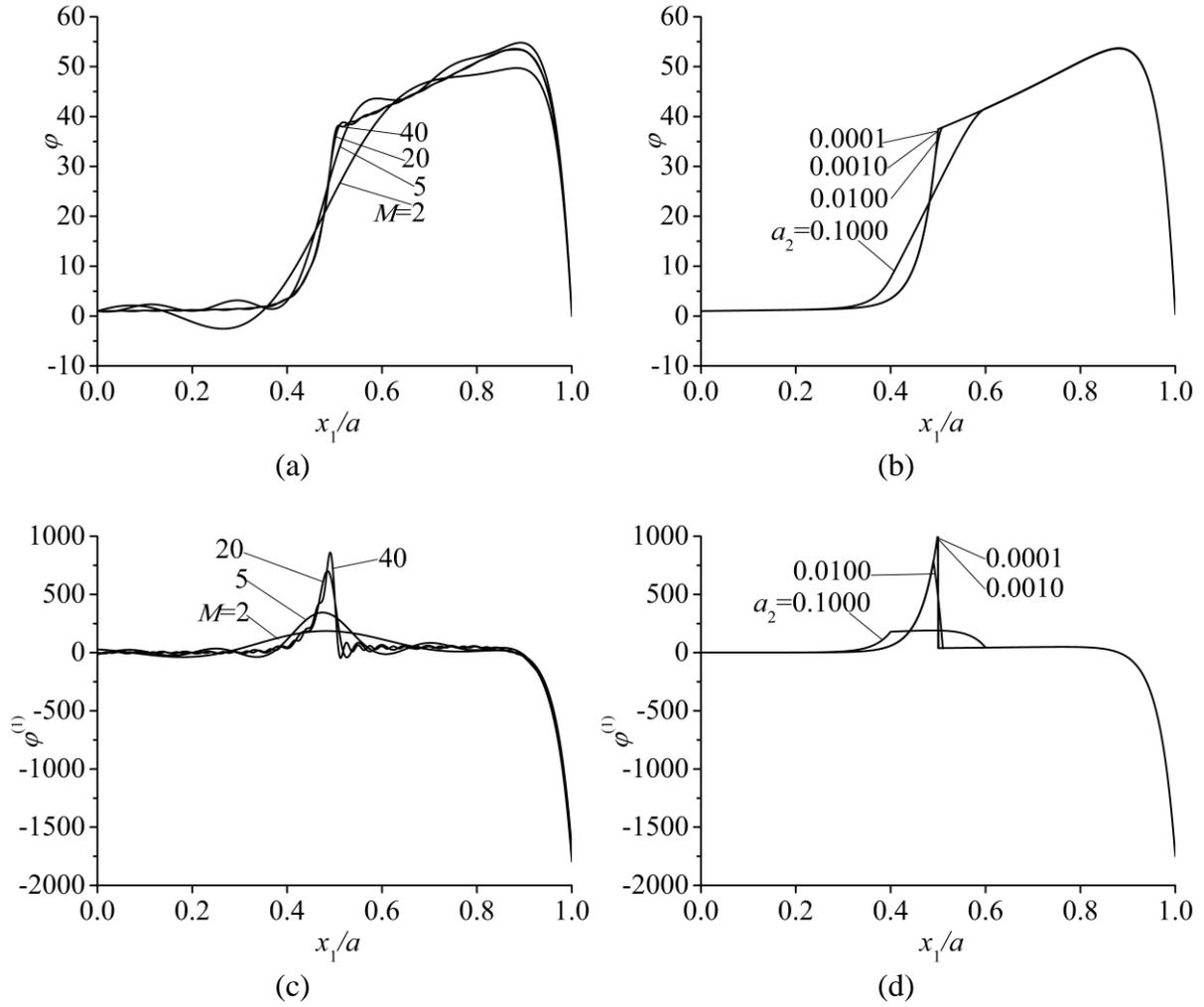

Figure 9: Two Convergence processes for the whole interval solution and the subinterval solution with computational parameters $P_e = 30$ and $D_a = 1$:

(a) $\varphi$ curves for the whole interval solution with different terms of $M$,

(b) $\varphi$ curves for the subinterval solution with different values of $a_2$,

(c) $\varphi^{(1)}$ curves for the whole interval solution with different terms of $M$,

(d) $\varphi^{(1)}$ curves for the subinterval solution with different values of $a_2$.

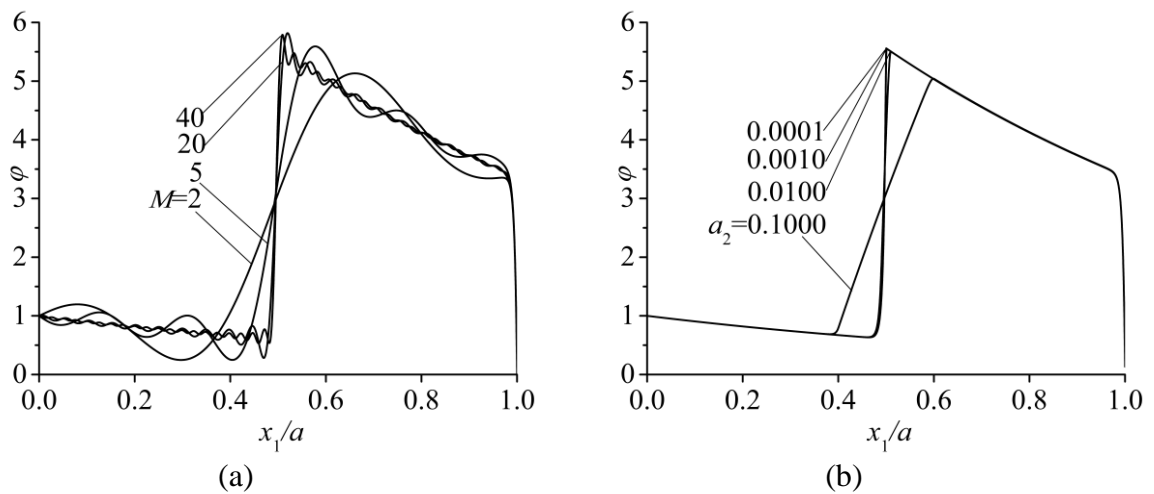

(a)                                             (b)

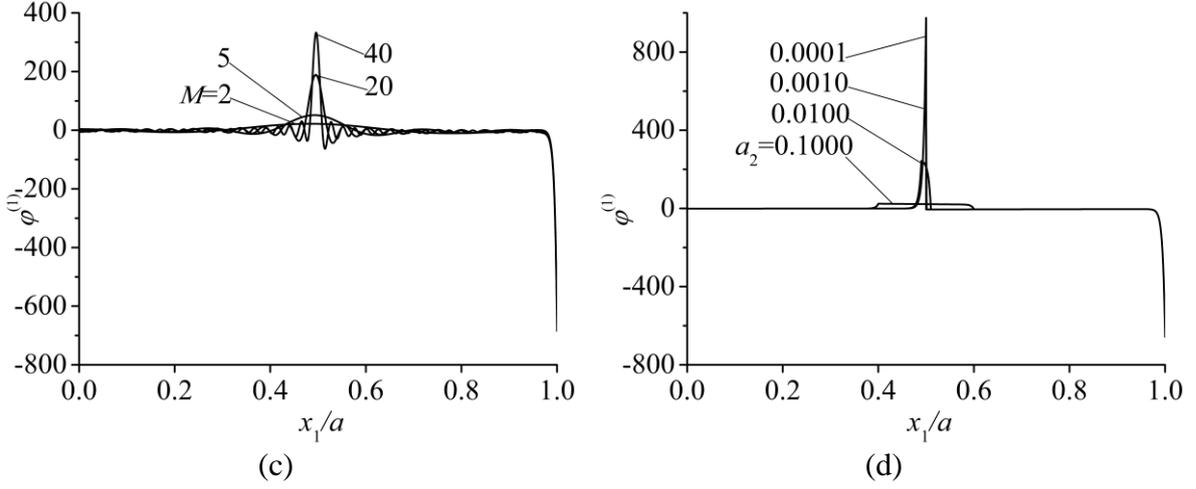

Figure 10: Two convergence processes for the whole interval solution and the subinterval solution with computational parameters $P_e = 200$ and $D_a = -1$:

(a) $\varphi$ curves for the whole interval solution with different terms of $M$,

(b) $\varphi$ curves for the subinterval solution with different values of $a_2$,

(c) $\varphi^{(1)}$ curves for the whole interval solution with different terms of $M$,

(d) $\varphi^{(1)}$ curves for the subinterval solution with different values of $a_2$.

## 4. Fourier series multiscale solution for two-dimensional convection-diffusion-reaction equation

The two-dimensional convection-diffusion-reaction equation is a second order ($2r = 2$) linear differential equation with constant coefficients, where convection refers to the term involving first order partial derivative of the unknown function. Therefore, we rationally specify the type of the Fourier series expansion of the two-dimensional solution as the full-range Fourier series and closely follow the step outlined in [39] to derive the Fourier series multiscale solution of the equation. For brevity, we are to offer a brief presentation of some necessary ingredients such as description of the problem, homogeneous solutions of the equation and corner function in the particular solution, rather than to give a systematic development of the Fourier series multiscale solution that consists of the general solution, the particular solution, and sometimes the supplementary solution.

1. Description of the problem

With reference to [20, 21, 41], the two-dimensional convection-diffusion-reaction equation in the domain $[-a, a] \times [-b, b]$ can be written in the dimensionless form as

$$\mathcal{L}_{cdr,2}\varphi = f, \tag{62}$$

where the differential operator

$$\mathcal{L}_{cdr,2} = P_{e1}\frac{\partial}{\partial x_1} + P_{e2}\frac{\partial}{\partial x_2} - (\frac{\partial^2}{\partial x_1^2} + \frac{\partial^2}{\partial x_2^2}) - P_e D_a, \tag{63}$$

$\varphi(x_1, x_2)$ is the solution function, $f(x_1, x_2)$ is a given source function, $P_e$ is the dimensionless Peclet number, $P_{e1}$ and $P_{e2}$ are the $x_1$ and $x_2$ components of $P_e$ respectively and depend on the inflow angle $\theta$, and $D_a$ is the dimensionless Damkohler number.

2. Homogeneous solutions

Suppose that Eq. (62) has homogeneous solutions of the following form

$$p_{1n,H}(x_1, x_2) = \exp(\eta_n x_1) \exp(\mathrm{i}\beta_n x_2) \quad (\mathrm{i} = \sqrt{-1}), \tag{64}$$

where $n$ is a nonnegative integer, $\eta_n$ is an undetermined constant, $\beta_n = n\pi/b$.

Substituting Eq. (64) into the homogeneous form of Eq. (62), we obtain the characteristic equation

$$\eta_n^2 - P_{e1}\eta_n + (-\beta_n^2 + P_e D_a - \mathrm{i}P_{e2}\beta_n) = 0. \tag{65}$$

If we write

$$\gamma_{1n} = -P_{e1}^2 - 4\beta_n^2 + 4P_e D_a, \quad \gamma_{2n} = -4P_{e2}\beta_n,$$

$$\alpha_{1n} = -\frac{P_{e1}}{2}, \quad \alpha_{2n} + \mathrm{i}\alpha_{3n} = \frac{1}{2}\sqrt{\gamma_{1n} + \mathrm{i}\gamma_{2n}}, \tag{66}$$

then we further formulate the characteristic roots as

$$\eta_{n,1} = (-\alpha_{1n} + \alpha_{3n}) - \mathrm{i}\alpha_{2n} \text{ and } \eta_{n,2} = (-\alpha_{1n} - \alpha_{3n}) + \mathrm{i}\alpha_{2n}. \tag{67}$$

It is obvious that, for $n = 0$, Eq. (65) has one double real root when $\gamma_{1n} = 0$, two distinct complex roots when $\gamma_{1n} > 0$, two distinct real roots when $\gamma_{1n} < 0$; for $n > 0$, Eq. (65) has one double real root when $\gamma_{1n} = 0$ and $\gamma_{2n} = 0$, and two distinct complex roots when $\gamma_{1n} \neq 0$ or $\gamma_{2n} \neq 0$.

Accordingly, for $n = 0$ Table 7 presents the sequence of homogeneous solutions $\{p_{10l,H}(x_1, x_2)\}_{1 \leq l \leq 2}$, and for $n > 0$ Table 8 presents the sequence of homogeneous solutions $\{p_{1nl,H}(x_1, x_2)\}_{1 \leq l \leq 4}$.

Table 7: Expressions for $\{p_{10l,H}(x_1, x_2)\}_{1 \leq l \leq 2}$.

| | $\gamma_{1n} = 0$ | $\gamma_{1n} > 0$ | $\gamma_{1n} < 0$ |
|---|---|---|---|
| $p_{101,H}(x_1, x_2)$ | $\dfrac{\exp(-\alpha_{1n}x_1)}{\cosh(-\alpha_{1n}a)}$ | $\dfrac{\exp(-\alpha_{1n}x_1)}{\cosh(-\alpha_{1n}a)}\cos(\alpha_{2n}x_1)$ | $\dfrac{\exp[(-\alpha_{1n} + \alpha_{3n})x_1]}{\cosh[(-\alpha_{1n} + \alpha_{3n})a]}$ |
| $p_{102,H}(x_1, x_2)$ | $\dfrac{x_1}{a}\dfrac{\exp(-\alpha_{1n}x_1)}{\cosh(-\alpha_{1n}a)}$ | $\dfrac{\exp(-\alpha_{1n}x_1)}{\cosh(-\alpha_{1n}a)}\sin(\alpha_{2n}x_1)$ | $\dfrac{\exp[(-\alpha_{1n} - \alpha_{3n})x_1]}{\cosh[(-\alpha_{1n} - \alpha_{3n})a]}$ |

Table 8: Expressions for $\{p_{1nl,H}(x_1, x_2)\}_{1 \leq l \leq 4}$.

| | $\gamma_{1n} = 0$ and $\gamma_{2n} = 0$ | $\gamma_{1n} \neq 0$ or $\gamma_{2n} \neq 0$ |
|---|---|---|
| $p_{1n1,H}(x_1, x_2)$ | $\dfrac{\exp(-\alpha_{1n}x_1)}{\cosh(-\alpha_{1n}a)}\cos(\beta_n x_2)$ | $\dfrac{\exp[(-\alpha_{1n} + \alpha_{3n})x_1]}{\cosh[(-\alpha_{1n} + \alpha_{3n})a]}\cos(\beta_n x_2 - \alpha_{2n}x_1)$ |
| $p_{1n2,H}(x_1, x_2)$ | $\dfrac{x_1}{a}\dfrac{\exp(-\alpha_{1n}x_1)}{\cosh(-\alpha_{1n}a)}\cos(\beta_n x_2)$ | $\dfrac{\exp[(-\alpha_{1n} + \alpha_{3n})x_1]}{\cosh[(-\alpha_{1n} + \alpha_{3n})a]}\sin(\beta_n x_2 - \alpha_{2n}x_1)$ |
| $p_{1n3,H}(x_1, x_2)$ | $\dfrac{\exp(-\alpha_{1n}x_1)}{\cosh(-\alpha_{1n}a)}\sin(\beta_n x_2)$ | $\dfrac{\exp[(-\alpha_{1n} - \alpha_{3n})x_1]}{\cosh[(-\alpha_{1n} - \alpha_{3n})a]}\cos(\beta_n x_2 + \alpha_{2n}x_1)$ |
| $p_{1n4,H}(x_1, x_2)$ | $\dfrac{x_1}{a}\dfrac{\exp(-\alpha_{1n}x_1)}{\cosh(-\alpha_{1n}a)}\sin(\beta_n x_2)$ | $\dfrac{\exp[(-\alpha_{1n} - \alpha_{3n})x_1]}{\cosh[(-\alpha_{1n} - \alpha_{3n})a]}\sin(\beta_n x_2 + \alpha_{2n}x_1)$ |

3. Structural characteristics of corner function in the particular solution

It is worth noting that Eq. (62) is a second order differential equation, and therefore the vector of basis functions $\mathbf{\Phi}_3^T$ and the vector of undetermined coefficients $\mathbf{q}_3$ actually consist of a single element, namely

$$\mathbf{\Phi}_3^T = [\varphi_{3,00}(x_1, x_2)], \tag{68}$$

$$\mathbf{q}_3^T = [q_{3,00}], \tag{69}$$

where the function

$$\varphi_{3,00}(x_1, x_2) = \frac{x_1 x_2}{4ab}, \tag{70}$$

and $q_{3,00}$ is an undetermined constant.

## 5. Two-dimensional numerical examples

In this section, we are to adopt an inverse validation method to the numerical analysis of the two-dimensional convection-diffusion-reaction equation. For this propose, we begin with an appropriate selection of a reference solution $\varphi_{ref}(x_1, x_2)$, and by plugging it back into Eq. (62), we implement an imposition of the source function and the boundary condition on the two-dimensional convection-diffusion-reaction equation. Based on this arrangement, we solve the two-dimensional convection-diffusion-reaction equation with the computational procedure outlined in section 4 and validate the obtained Fourier series multiscale solution through a comparison with the reference solution.

For simplicity, we take the following linear combination of homogeneous solutions $\{p_{1nl,H}(x_1, x_2)\}_{1 \leq l \leq 4}$ as the reference solution

$$\varphi_{ref}(x_1, x_2) = \sum_{l=1}^{4} p_{1nl,H}(x_1, x_2), \tag{71}$$

where the parameter $\beta_n = n\pi/b$ is substituted with $\beta_{ref} = \pi/2b$.

And accordingly, we find that the source function

$$f_{ref}(x_1, x_2) = 0. \tag{72}$$

It is also worthy of notice that, in the settings of computational scheme of the Fourier series multiscale solution of the two-dimensional convection-diffusion-reaction equation, the order of the interpolation algebraical polynomial of the source function is zero and the derivation method for discrete equations is the Fourier coefficient comparison method.

*5.1. Convergence characteristics*

As shown in Table 9, we perform four comparative convergence experiments in series for detailed analysis of the influences of some possible factors of the two-dimensional convection-diffusion-reaction equation, such as the computational parameters, boundary condition, inflow angle and length-width ratio, on convergence characteristics of the Fourier series multiscale solution. For instance, in the first comparative convergence experiment, we start out with the reference computational scheme (see the first column of Table 9 ), and change the computational parameters from strong reaction type to reaction-dominated type, convection-dominated type and strong convection type successively; in the second comparative convergence experiment, we start out with the reference computational scheme, and change the inflow angle from $\pi/3$ to $\pi/4$, $\pi/6$ and zero successively; in the third comparative convergence experiment, we start out with the reference computational scheme, and change the length-width ratio from 1.0 to 0.67, 0.50, 1.25 and 2.0 successively; and in the fourth comparative convergence experiment, we start out with the reference computational scheme, and change the boundary condition from the DDDD boundary condition to the DDND boundary condition and the DNND boundary condition successively (see Figure 11).

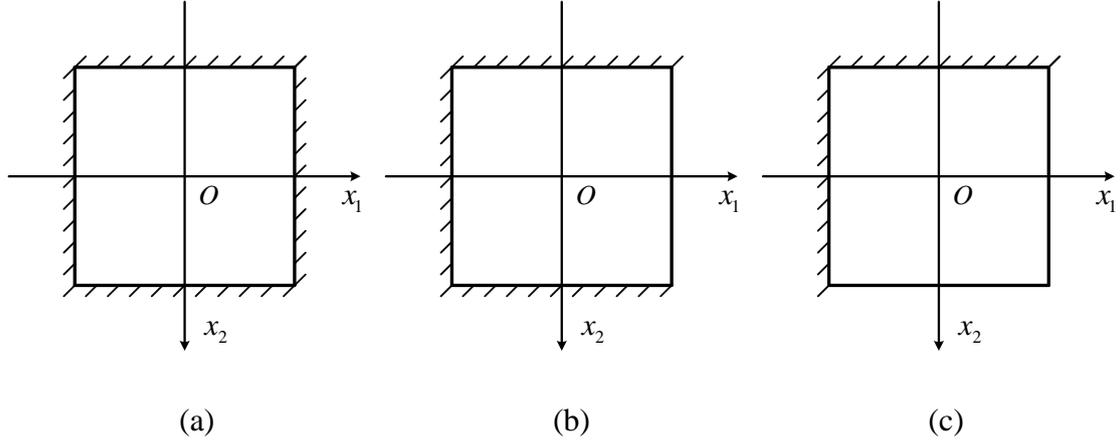

(a)            (b)            (c)

Figure 11: The types of boundary conditions considered: (a) DDDD boundary condition, (b) DDND boundary condition, (c) DNND boundary condition.

Table 9: Comparative convergence experiments for two-dimensional convection-diffusion-reaction equation.

| Comparative convergence experiment | No. | Boundary condition | Computational parameter $(P_e, D_a)$ | Inflow angle $\theta$ | length width ratio $a/b$ |
|---|---|---|---|---|---|
| 1 | a | DDDD | (3, 90) | $\pi/3$ | 1.0 |
|   | b |      | (1, 30) |         |     |
|   | c |      | (30, 1) |         |     |
|   | d |      | (200, -1) |       |     |
| 2 | a | DDDD | (3, 90) | $\pi/3$ | 1.0 |
|   | b |      |         | $\pi/4$ |     |
|   | c |      |         | $\pi/6$ |     |
|   | d |      |         | 0       |     |
| 3 | a | DDDD | (3, 90) | $\pi/3$ | 1.0 |
|   | b |      |         |         | 0.67 |
|   | c |      |         |         | 0.50 |

|   |   |   |   |   | 1.25 |
|---|---|---|---|---|------|
|   | d |   |   |   |      |
|   | e |   |   |   | 2.0  |
|   | a | DDDD |   |   |    |
| 4 | b | DDND | (3, 90) | $\pi/3$ | 1.0 |
|   | c | DNND |   |   |    |

As to the reference computational scheme and the four comparative convergence experiments given in Table 9, we truncate the composite Fourier series of the Fourier series multiscale solution successively with the first 2, 3, 5, 10, 20, 30 and 40 terms and then we compute the overall computational errors, internal computational errors, boundary computational errors and corner computational errors of the function $\varphi(x_1, x_2)$ and its first order partial derivatives. Some of the results are shown in Figures 12-17. We make a brief analysis of the convergence characteristics of the Fourier series multiscale solution of the two-dimensional convection-diffusion-reaction equation as the following:

a. As shown in Figure 12, the Fourier series multiscale solution has good convergence in general. And especially, the composite Fourier series of $\varphi(x_1, x_2)$ and its first order partial derivatives converge well within the solution domain, the composite Fourier series of $\varphi(x_1, x_2)$ and its first order partial derivative $\varphi^{(1,0)}(x_1, x_2)$ converge well on the boundary of the solution domain, the composite Fourier series of $\varphi^{(0,1)}(x_1, x_2)$, the first order partial derivative of $\varphi(x_1, x_2)$, does not converge as well on the boundary of the solution domain, the composite Fourier series of $\varphi(x_1, x_2)$ converge well at the corners of the solution domain, and the composite Fourier series of the first order partial derivatives of $\varphi(x_1, x_2)$ do not converge as well at the corners of the solution domain.

b. As shown in Figures 13-16, the influence factors have different effects on the convergence characteristics of the Fourier series multiscale solution. For instance, the adjustment of the computational parameters has marked influences on the convergence characteristics of the Fourier series multiscale solution, the adjustment of the inflow angle bring about few effects on the convergence of the Fourier series multiscale solution, the decrease of length-width ratio significantly worsens the convergence of the Fourier series multiscale solution, the increase of length-width ratio has few effects on the convergence of the Fourier series multiscale solution, and the adjustment of the boundary condition from the DDDD boundary condition to the DDND boundary condition or the DNND boundary condition markedly compromises the convergence of the Fourier series multiscale solution.

c. As shown in Figure 17, with the length-width ratio specified as $a/b = 0.50$, the adjustment of $N/M$, the ratio of number of truncated terms of the composite Fourier series, from 1.0 to 2.0 significantly improves the convergence characteristics of Fourier series multiscale solution.

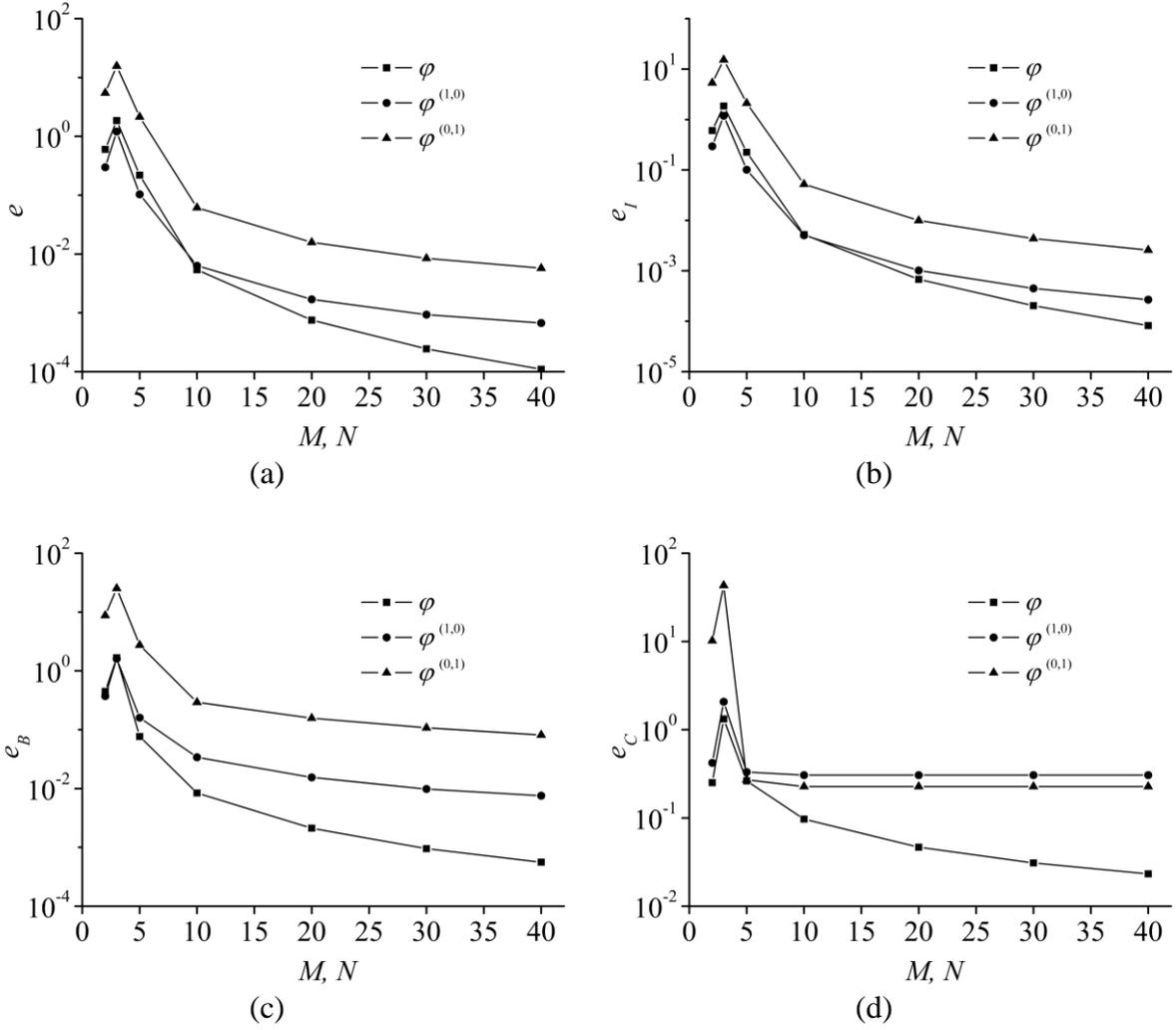

Figure 12: Convergence characteristics of the Fourier series multiscale solution of two-dimensional convection-diffusion-reaction equation (the reference configuration):
(a) $e^{(k_1,k_2)}(\varphi_{M,N})$ - $M$, $N$ curves, (b) $e_I^{(k_1,k_2)}(\varphi_{M,N})$ - $M$, $N$ curves, (c) $e_B^{(k_1,k_2)}(\varphi_{M,N})$ - $M$, $N$ curves, (d) $e_C^{(k_1,k_2)}(\varphi_{M,N})$ - $M$, $N$ curves.

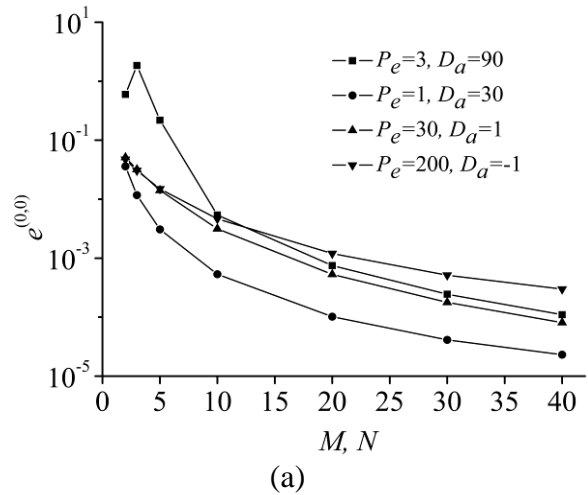

(a)

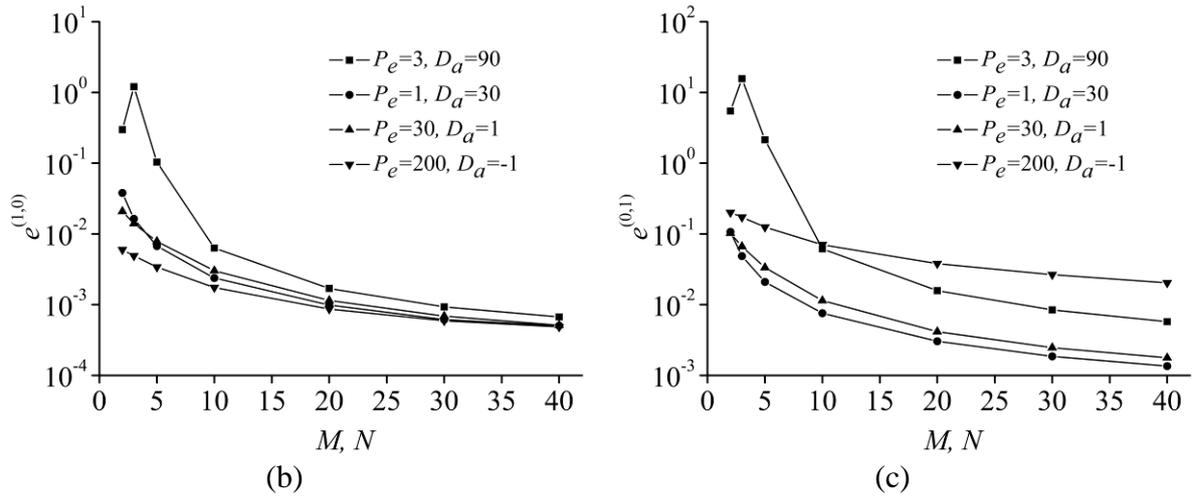

(b)                      (c)

Figure 13: Convergence comparison of the Fourier series multiscale solutions with different computational parameters $P_e$ and $D_a$:
(a) $e^{(0,0)}(\varphi_{M,N})$ - $M$, $N$ curves, (b) $e^{(1,0)}(\varphi_{M,N})$ - $M$, $N$ curves, (c) $e^{(0,1)}(\varphi_{M,N})$ - $M$, $N$ curves.

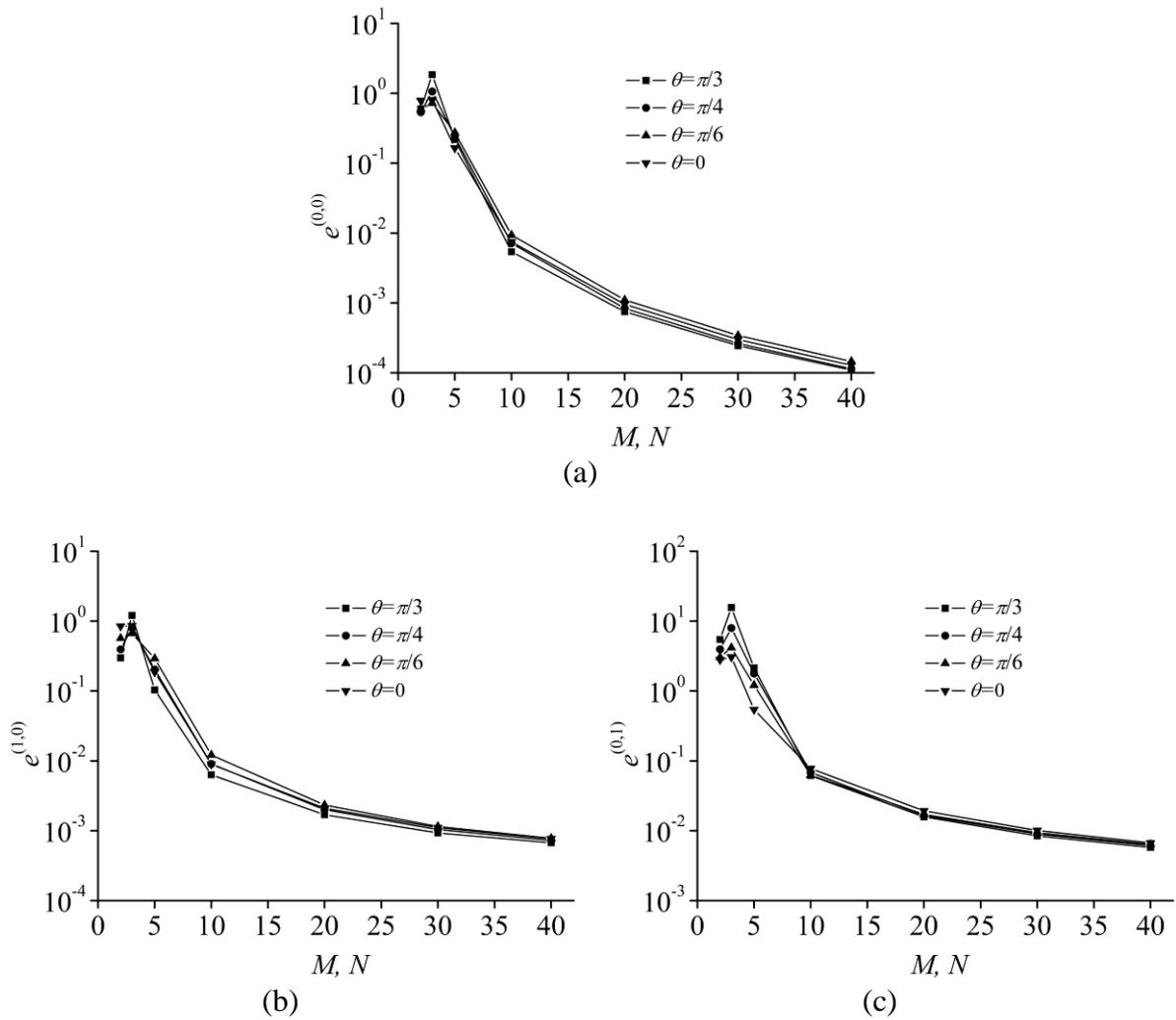

(b)                      (c)

Figure 14: Convergence comparison of the Fourier series multiscale solutions with different inflow angles:

(a) $e^{(0,0)}(\varphi_{M,N})$ - $M$, $N$ curves, (b) $e^{(1,0)}(\varphi_{M,N})$ - $M$, $N$ curves, (c) $e^{(0,1)}(\varphi_{M,N})$ - $M$, $N$ curves.

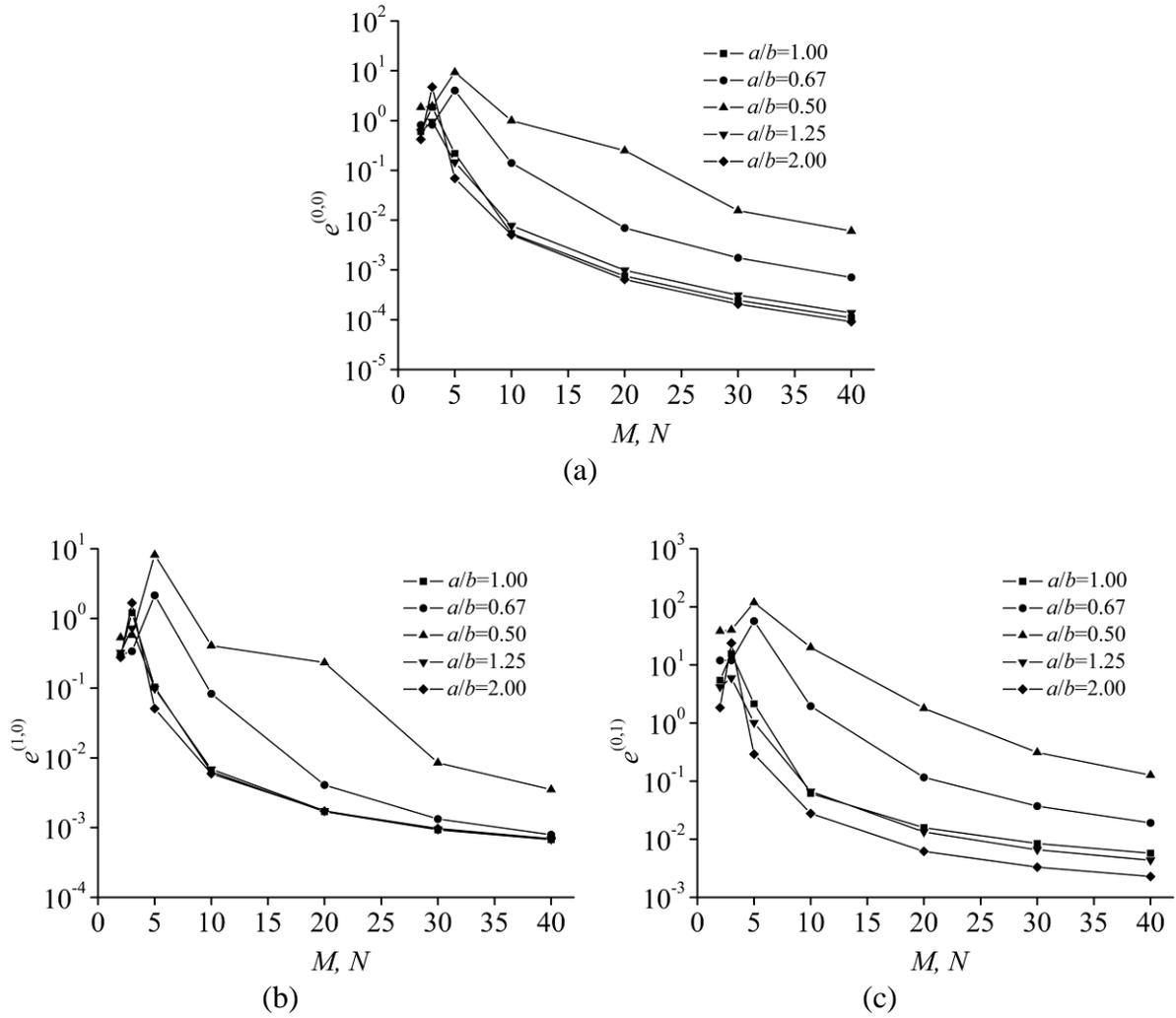

Figure 15: Convergence comparison of the Fourier series multiscale solutions with different length-width ratios:
(a) $e^{(0,0)}(\varphi_{M,N})$ - $M$, $N$ curves, (b) $e^{(1,0)}(\varphi_{M,N})$ - $M$, $N$ curves, (c) $e^{(0,1)}(\varphi_{M,N})$ - $M$, $N$ curves.

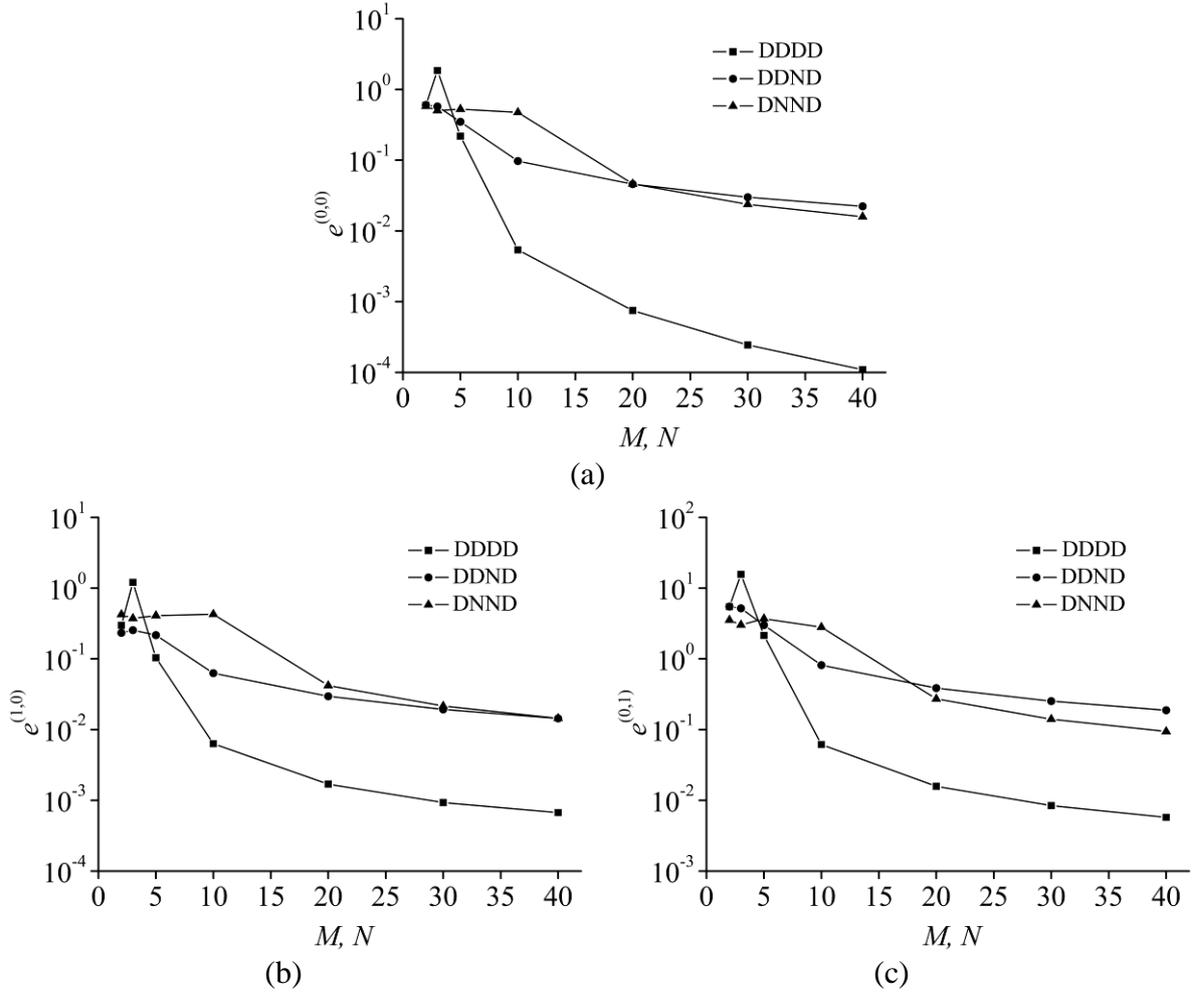

Figure 16: Convergence comparison of the Fourier series multiscale solutions with different boundary conditions:
(a) $e^{(0,0)}(\varphi_{M,N})$ - $M$, $N$ curves, (b) $e^{(1,0)}(\varphi_{M,N})$ - $M$, $N$ curves, (c) $e^{(0,1)}(\varphi_{M,N})$ - $M$, $N$ curves.

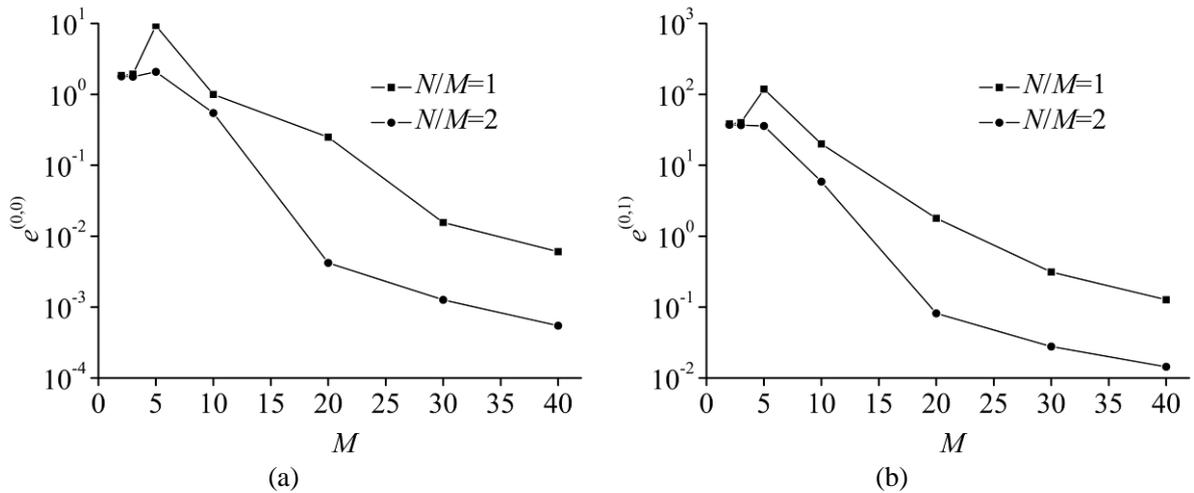

Figure 17: Convergence comparison of the Fourier series multiscale solutions with different values of $N/M$ (length-width ratio $a/b = 0.5$):
(a) $e^{(0,0)}(\varphi_{M,N})$ - $M$ curves, (b) $e^{(0,1)}(\varphi_{M,N})$ - $M$ curves.

*5.2. Multiscale characteristics*

As shown in Table 9, by starting with the reference computational scheme specified in the comparative convergence experiments and changing the computational parameters of the two-dimensional convection-diffusion-reaction equation from strong reaction type to reaction-dominated type, convection-dominated type and strong convection type successively, we present in Figures 18-21 the corresponding reference results and the computed results of the Fourier series multiscale method with the numbers of truncated terms of the composite Fourier series specified as $M = N = 40$. It is demonstrated that, with the adjustment of the computational parameters, the computed results obtained by the Fourier series multiscale method is identical with the reference results; the distribution patterns throughout the solution domain evolves from periodic oscillation to exponential decay; and when the computational parameters turn to be the strong convection type, the boundary layer appears as a typical multiscale phenomenon in the solution domain.

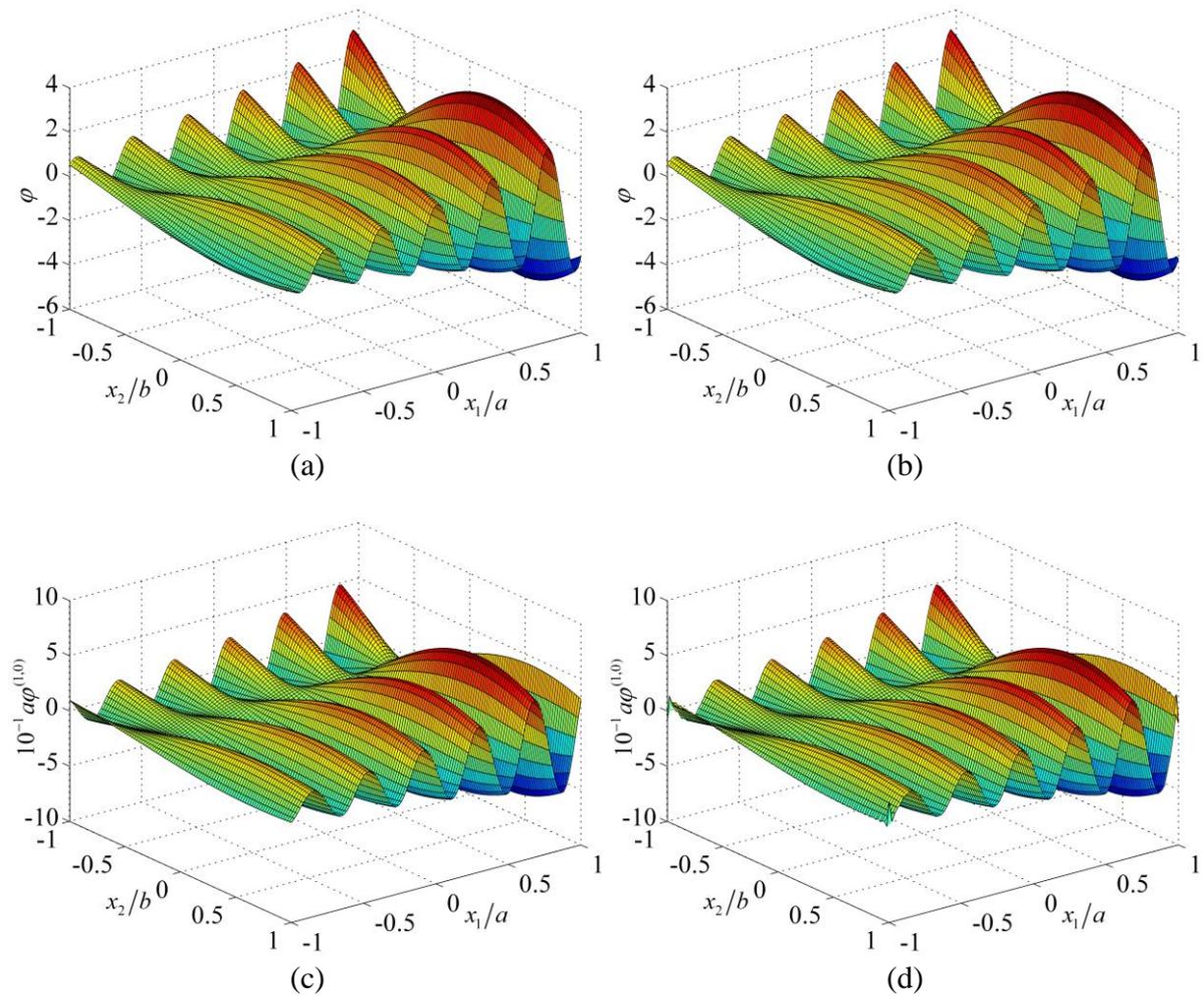

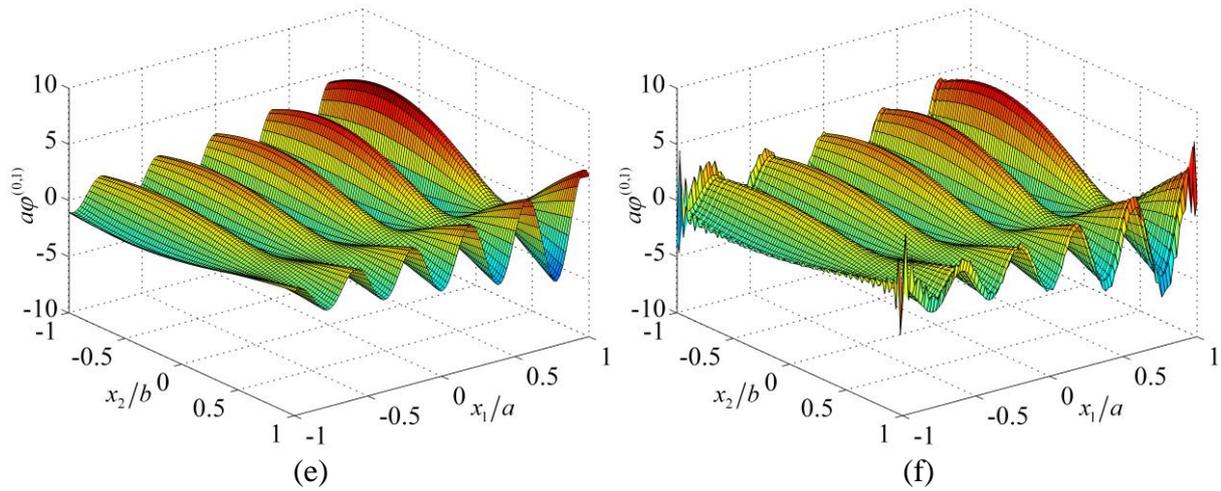

Figure 18: Results for the computational parameters $P_e = 3$ and $D_a = 90$:
(a) reference result of $\varphi(x_1, x_2)$, (b) computed result of $\varphi(x_1, x_2)$,
(c) reference result of $\varphi^{(1,0)}(x_1, x_2)$, (d) computed result of $\varphi^{(1,0)}(x_1, x_2)$,
(e) reference result of $\varphi^{(0,1)}(x_1, x_2)$, (f) computed result of $\varphi^{(0,1)}(x_1, x_2)$.

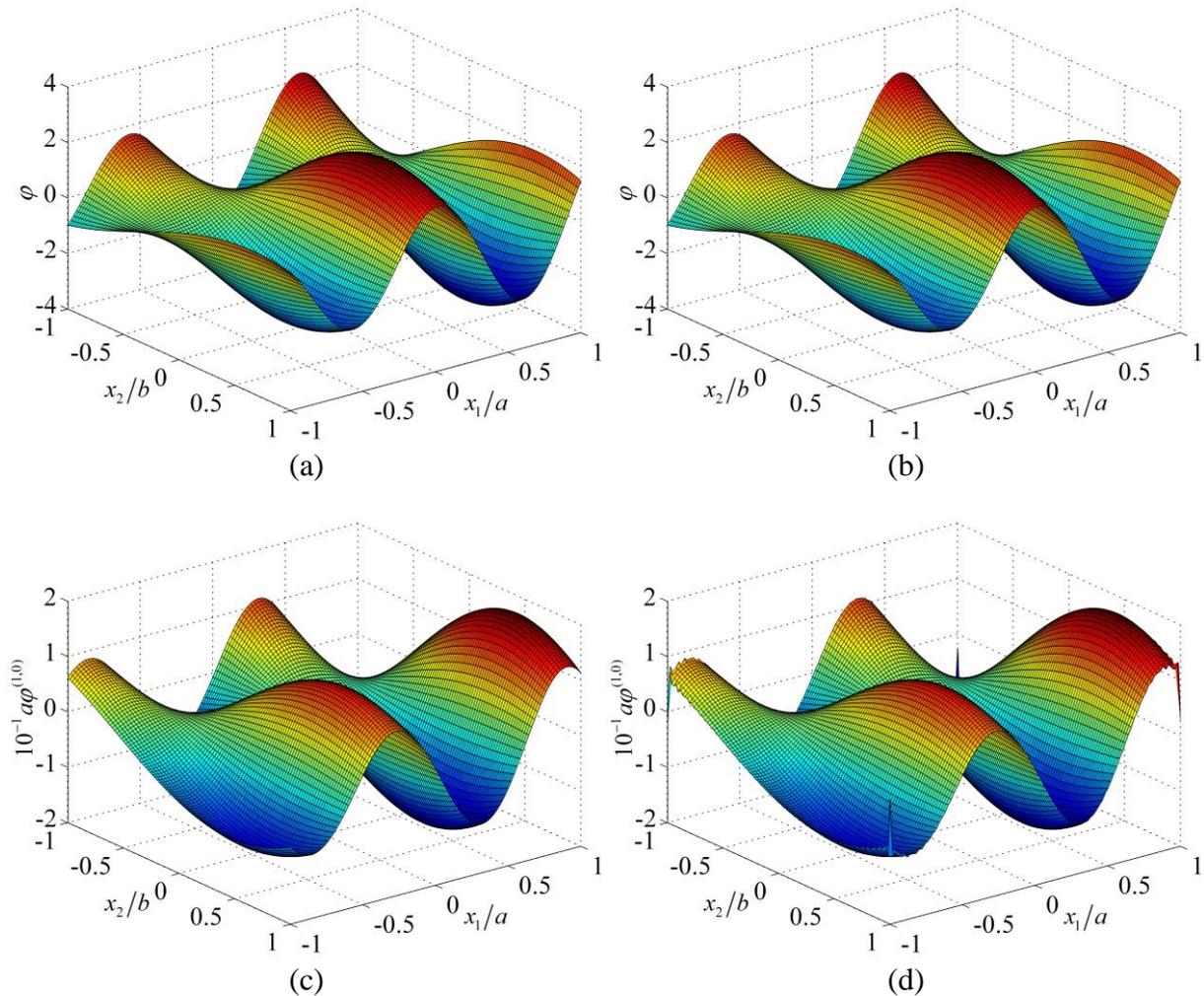

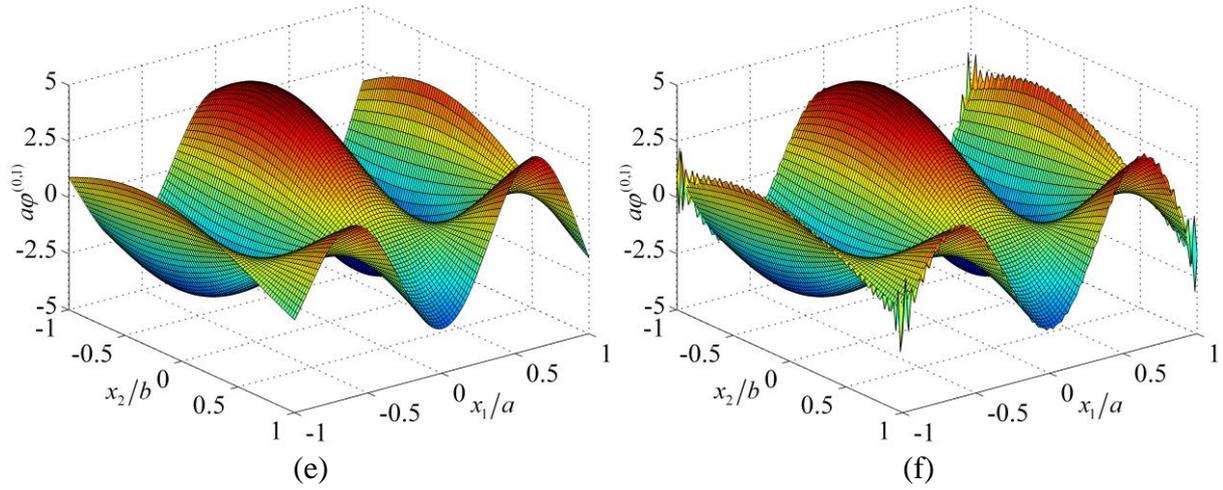

Figure 19: Results for the computational parameters $P_e = 1$ and $D_a = 30$:
(a) reference result of $\varphi(x_1, x_2)$, (b) computed result of $\varphi(x_1, x_2)$,
(c) reference result of $\varphi^{(1,0)}(x_1, x_2)$, (d) computed result of $\varphi^{(1,0)}(x_1, x_2)$,
(e) reference result of $\varphi^{(0,1)}(x_1, x_2)$, (f) computed result of $\varphi^{(0,1)}(x_1, x_2)$.

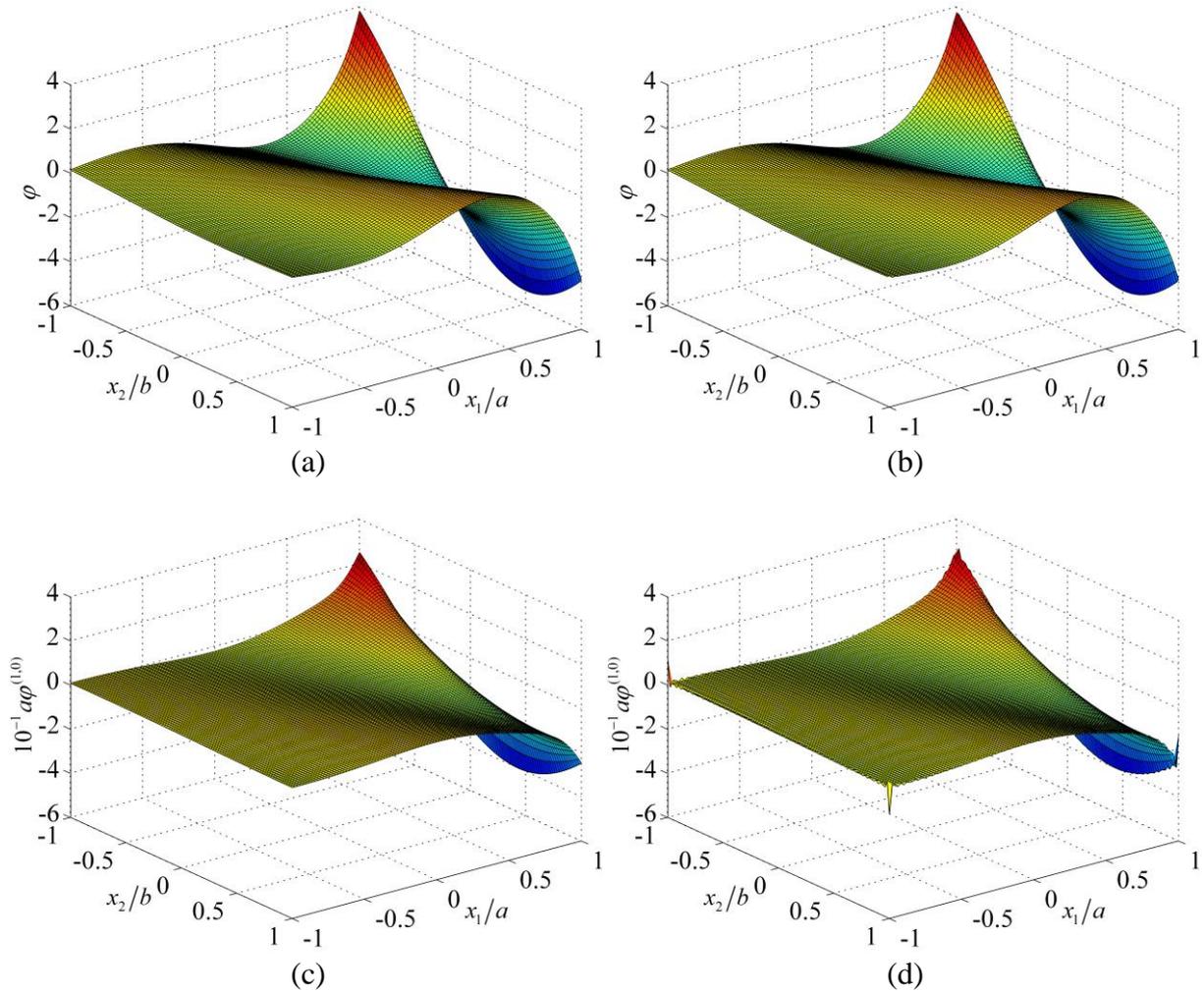

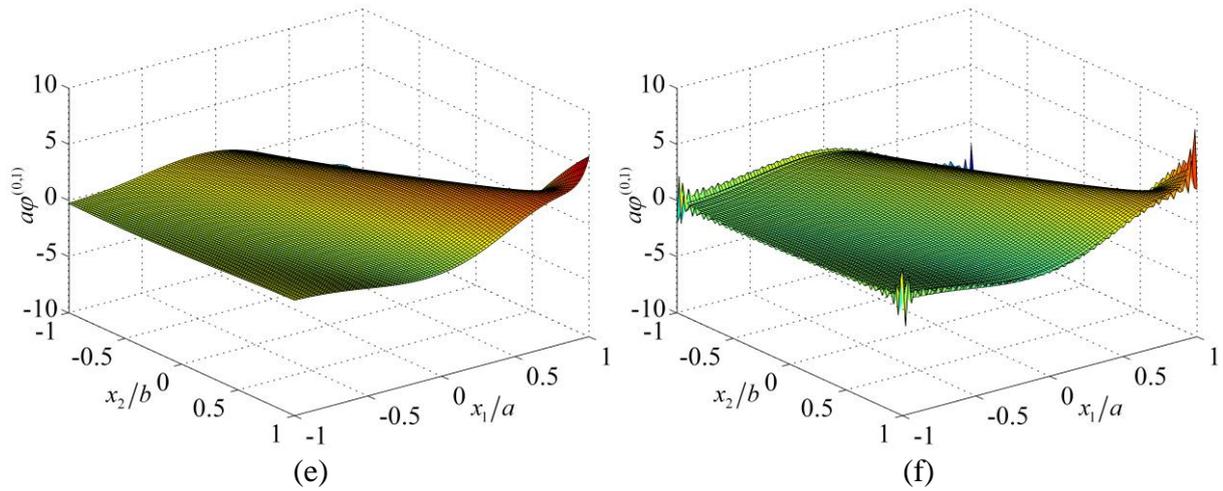

(e)                                         (f)

Figure 20: Results for the computational parameters $P_e = 30$ and $D_a = 1$:
(a) reference result of $\varphi(x_1, x_2)$, (b) computed result of $\varphi(x_1, x_2)$,
(c) reference result of $\varphi^{(1,0)}(x_1, x_2)$, (d) computed result of $\varphi^{(1,0)}(x_1, x_2)$,
(e) reference result of $\varphi^{(0,1)}(x_1, x_2)$, (f) computed result of $\varphi^{(0,1)}(x_1, x_2)$.

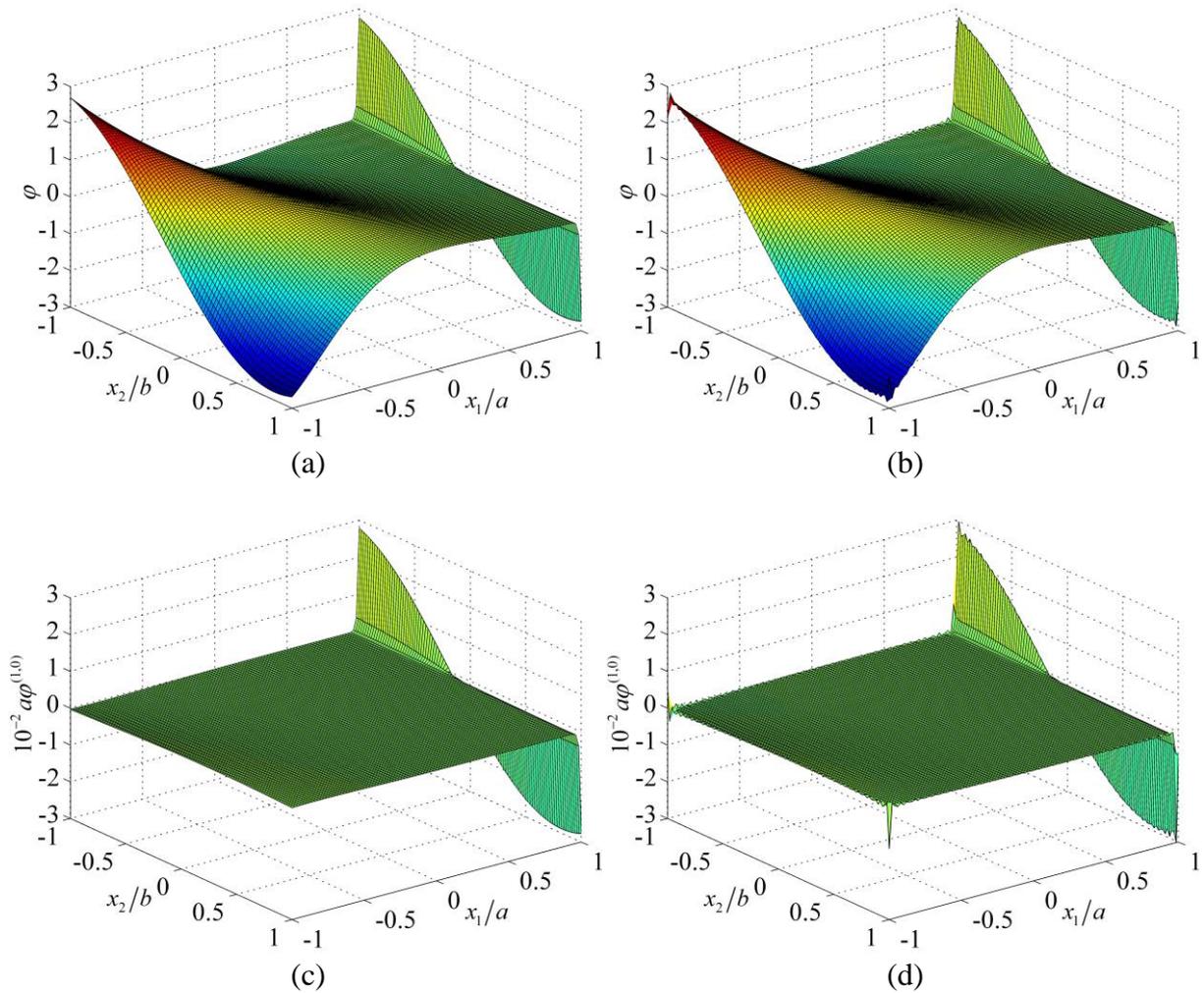

Figure 21: Results for the computational parameters $P_e = 200$ and $D_a = -1$:

(a) reference result of $\varphi(x_1, x_2)$, (b) computed result of $\varphi(x_1, x_2)$,
(c) reference result of $\varphi^{(1,0)}(x_1, x_2)$, (d) computed result of $\varphi^{(1,0)}(x_1, x_2)$,

## 6. Conclusions

The convection-diffusion-reaction equation is a benchmark problem for the multiscale methods. In this paper, we perform thorough investigation on the application of the Fourier series multiscale method to the one-dimensional and two-dimensional convection-diffusion-reaction equations. It is concluded that

1. We derive the Fourier series multiscale solutions of the one-dimensional and two-dimensional convection-diffusion-reaction equations, and investigate the convergence characteristics of the obtained Fourier series multiscale solutions.

2. We bring to an integration of the Fourier series multiscale method and the derivation technique for discrete equations in solution of the one-dimensional and two-dimensional convection-diffusion-reaction equations, and demonstrate the merit of the Fourier series multiscale method which yields stable and accurate numerical results for wide range of computational parameters and boundary conditions.

3. We reveal the multiscale characteristics of the one-dimensional and two-dimensional convection-diffusion-reaction equations.

The preliminary study on applications verifies the effectiveness of the present Fourier series multiscale method and provides a reliable reference which can be used for persistent improvement in computational performance of other multiscale methods.